\documentclass[12pt]{article}
\usepackage{latexsym,amsmath,amssymb}

\newif\ifdviwin

\usepackage[latin1]{inputenc}
\usepackage[english]{babel}
\usepackage{indentfirst}
\usepackage[mathscr]{eucal}
\usepackage{amssymb,amsmath,amsfonts}
\usepackage{graphicx}

\newif\ifdviwin

\dviwintrue

\def\cG{\mathcal{G}}

\def\cN{\mathcal{N}}

\def\vt{\vartheta}
\let\8=\infty \let\0=\emptyset  
\let\hat=\widehat

\let\tilde=\widetilde

\let\landa=\lambda
\let\alfa=\alpha

\let\parc=\partial

\def\ep{\varepsilon}
\def\landa{\lambda}

\def\flecha{\rightarrow}
\def\esiz{\langle}
\def\esde{\rangle}

\def\cte.{\mathop{\rm cte.}\nolimits}

\def\Re{\mathop{\rm Re }\nolimits}

\def\cosh{\mathop{\rm cosh }\nolimits}

\def\N{\mathbb{N}}
\def\L{\mathbb{L}}

\def\R{\mathbb{R}}
\def\Z{\mathbb{Z}}
\def\C{\mathbb{C}}
\def\D{\mathbb{D}}

\def\H{\mathbb{H}}
\def\S{\mathbb{S}}

 \newtheorem{defi}{Definition}
 \newtheorem{teo}[defi]{Theorem}
 \newtheorem{pro}[defi]{Proposition}
 \newtheorem{cor}[defi]{Corollary}
 \newtheorem{lem}[defi]{Lemma}

 \newenvironment{proof}{\rm \trivlist \item[\hskip \labelsep{\it
      Proof}:]}{\par\nopagebreak \hfill $\Box$ \endtrivlist}

 \newenvironment{proof1}{\rm \trivlist \item[\hskip \labelsep{\it
      Proof of Theorem \ref{th:main}}:]}{\par\nopagebreak \hfill $\Box$ \endtrivlist}

\numberwithin{equation}{section}

\begin{document}
\mbox{}\vspace{0.4cm}\mbox{}

\begin{center}
\rule{14cm}{1.5pt}\vspace{0.5cm}

{\Large \bf Harmonic maps and constant mean curvature} \\ [0.3cm]{\Large \bf surfaces in $\H^2 \times \R$}\\
\vspace{0.5cm} {\large Isabel Fernández$\mbox{}^a$ and Pablo Mira$\mbox{}^b$}\\ \vspace{0.3cm}
\rule{14cm}{1.5pt}
\end{center}
  \vspace{1cm}
$\mbox{}^a$ Departamento de Geometría y Topología, Universidad de Granada,
E-18071 Granada, Spain. \\ e-mail: isafer@ugr.es \vspace{0.2cm}

\noindent $\mbox{}^b$ Departamento de Matemática Aplicada y Estadística,
Universidad Politécnica de Cartagena, E-30203 Cartagena, Murcia, Spain. \\
e-mail: pablo.mira@upct.es \vspace{0.2cm}

\vspace{0.2cm}


\vspace{0.3cm}

 \begin{abstract}
We introduce a \emph{hyperbolic Gauss map} into the Poincaré disk for any surface in $\H^2\times \R$ with regular vertical projection, and prove that if the surface has constant mean curvature $H=1/2$, this hyperbolic Gauss map is harmonic. Conversely, we show that every nowhere conformal harmonic map from an open simply connected Riemann surface $\Sigma$ into the Poincaré disk is the hyperbolic Gauss map of a two-parameter family of such surfaces. As an application we obtain that any holomorphic quadratic differential on $\Sigma$ can be realized as the Abresch-Rosenberg holomorphic differential
of some, and generically infinitely many, complete surfaces with $H=1/2$ in $\H^2 \times \R$. A similar result applies to minimal surfaces in the Heisenberg group
${\rm Nil_3}$. Finally, we classify all complete minimal vertical graphs in $\H^2\times \R$.
 \end{abstract}

\section{Introduction}

It is a classical result that the Hopf differential of a surface of constant mean curvature (CMC) in $\R^3$, and more generally in any $3$-dimensional space form, is holomorphic. Regarding more general target spaces, Abresch and Rosenberg \cite{AbRo} have recently proved that, even though the usual Hopf differential fails to be holomorphic on CMC-surfaces  in $\H^2 \times \R$ and $\S^2 \times \R$, a certain \emph{perturbed} Hopf differential is always holomorphic on any surface of this type. This striking theorem has put the topic into a new light, and is inspiring many research works on CMC surfaces in general ambient spaces \cite{Abr,BeTa,Dan,HLR,NeRo1,NeRo2,SaE,SaTo} (see also \cite{ACT,AEG}).

On the other hand, a condition stronger than the holomorphicity of the Hopf differential holds on any CMC surface in $\R^3$: the Gauss map of a CMC surface in $\R^3$ is a harmonic map into $\S^2$. Let us recall that if $h:\Sigma\flecha (M^2,\esiz, \esde)$ is a harmonic map from a Riemann surface $\Sigma$ into a Riemannian surface $M^2$, then the quantity $\esiz h_z,h_z\esde dz^2$ is a holomorphic quadratic differential on $\Sigma$. So any harmonic map into $\S^2$ comes along with a holomorphic differential, and in the case of CMC surfaces in $\R^3$ the differential induced by its harmonic Gauss map agrees up to a constant with the Hopf differential of the surface.

As a general rule of thought, the appearance of a geometrically defined harmonic map in the study of some class of surfaces should not be ignored, because it may allow to use techniques from the highly developed theory of integrable systems in the description of such surfaces. So, it is natural to analyze if the Abresch-Rosenberg holomorphic quadratic differential comes from some harmonic map geometrically defined on any CMC surface in $\H^2\times \R$ and $\S^2\times\R$.

The starting point of the present work is that for the special value $H=1/2$ of the mean curvature $H$ of a surface in $\H^2 \times \R$, a geometrically defined harmonic Gauss map into the Poincaré disk can be constructed. Furthermore the holomorphic quadratic differential associated to this harmonic map coincides up to a sign with the Abresch-Rosenberg holomorphic differential.

In the previous works on CMC surfaces in $\H^2 \times \R$ it has become clear that the  $H=1/2$ class is a limit case between two different situations. For instance, an embedded CMC surface in $\H^2 \times \R$ can be compact only if $H>1/2$ (\cite{HsHs,NeRo1}). In this sense, the mean curvature one half surfaces in $\H^2 \times \R$ are analogous to the surfaces with $H=1$ in $\H^3$, usually called \emph{Bryant surfaces}. These Bryant surfaces have a quite explicit form in terms of holomorphic data that is not shared by general CMC surfaces in $\H^3$, and their geometry has been studied in great detail (see for instance \cite{Bry,UY1,UY2,UY3,CHR,HRR,GaMi} and references therein). Another similarity between mean curvature one half surfaces in $\H^2 \times \R$ and Bryant surfaces is described in \cite{Dan} in terms of a Lawson-type correspondence for CMC surfaces in homogeneous spaces.

In this line, our results seem to confirm to some extent these analogies, as they suggest that the surfaces with $H=1/2$ in $\H^2 \times \R$ admit a more explicit treatment than the CMC surfaces with $H\neq 0, 1/2$ in $\H^2 \times \R$.

Our main purpose in this paper is to use the above harmonic Gauss map in order to study the surfaces with mean curvature one half in $\H^2 \times \R$. Our working scheme here goes along the lines of \cite{AbRo}. There, a geometric deep result (the classification of CMC spheres in product spaces) is obtained as a consequence of a known theorem from an independent theory (non-existence of non-zero holomorphic quadratic differentials on the sphere), once an adequate geometrical object (a Hopf type differential) has been found. In this work, we will use our construction of a harmonic Gauss map for surfaces with $H=1/2$ in $\H^2\times \R$ to translate into solutions to difficult geometric problems several known results from the theory of harmonic maps into the Poincaré disk.

We will also deal in this paper with minimal surfaces in $\H^2\times \R$. It is known that the horizontal factor $N$ as well as the height function $h$ of a minimal surface $\psi =(N,h)$ in $\H^2\times \R$ are harmonic maps into the Poincaré disk and the real line, respectively. So, this class fits into our discussion.

The geometry of minimal surfaces in spaces of the form $M^2\times \R$ has received a big number of contributions in recent years \cite{Ros,MeRo1,MeRo2,NeRo3,Hau,AlDa,SaTo,SaE,MMP}. The most studied situation is when $M^2$ has non-negative curvature, since in that case many clean classification theorems can be obtained. However, when the curvature of $M^2$ is negative, the situation changes radically. For instance, while complete minimal vertical graphs in $M^2\times \R$ are totally geodesic if $M$ has non-negative curvature \cite{Ros,AlDa}, there are many complete minimal vertical graphs in $\H^2\times \R$, as shown in \cite{NeRo3}. Our main result regarding minimal surfaces is the description of the space of all complete minimal vertical graphs in $\H^2\times \R$. Again, we will use for that known results from the theory of harmonic maps into the Poincaré disk.

Let us describe briefly our main results. Let $\psi =(N,h):\Sigma\flecha \H^2 \times \R\subset \L^4$ be a conformal immersion from a Riemann surface $\Sigma$ into $\H^2 \times \R$. Assume that $N:\Sigma\flecha \H^2$ is regular everywhere, and let $\eta =(\hat{N},u)$ be the unit normal to $\psi$. Then $u\neq 0$ at every point, and we can thus consider the map $\xi := (\eta +N)/u$, taking values in the intersection of the light cone $\N^3 \subset \L^4$ with the horizontal affine hyperplane $x_3=1$ of $\L^4$. Thus there is some $G:\Sigma\flecha \H^2$ such that $\xi =(G,1)$. We call $G$ the \emph{hyperbolic Gauss map} of $\psi$, because of its similarity with the usual hyperbolic Gauss map for surfaces in $\H^3$.

With the above construction, we prove in Section 3 that \emph{$G$ is a harmonic map on every mean curvature one half surface in $\H^2 \times \R$ with regular vertical projection}. Moreover, we shall show that in this case $\esiz G_z, G_z\esde dz^2 $ agrees with the Abresch-Rosenberg holomorphic differential up to a sign.

Also in Section 3 we deal with the inverse problem: is every harmonic map into $ \H^2$ the hyperbolic Gauss map of some mean curvature one half surface in $\H^2 \times \R$? We shall provide an answer to this question in the simply connected case, by characterizing the harmonic maps that are realized as hyperbolic Gauss maps. In particular, every nowhere conformal harmonic map into $\H^2$ belongs to this class. Furthermore, we will show that the class of surfaces with $H=1/2$ that share the hyperbolic Gauss map is a continuous $2$-parameter family, and we will determine when two such surfaces are congruent in $\H^2\times \R$. This property makes a big difference with the case of Bryant surfaces, where the space of surfaces with the same hyperbolic Gauss map can be infinite dimensional. We shall also provide an explicit formula that recovers from a harmonic map $G:\Sigma\flecha \H^2$ all the mean curvature one half surfaces in $\H^2 \times \R$ having $G$ as their hyperbolic Gauss map.

In Section 4 we will provide several applications of the previous results. The main one deals with the existence of complete examples, and is the following: \emph{any holomorphic quadratic differential on an open simply connected Riemann surface can be realized as the Abresch-Rosenberg differential of some, and generically infinitely many, complete surfaces with $H=1/2$ in $\H^2 \times \R$.} Prior to this result, the only known examples of complete non-minimal CMC surfaces in $\H^2 \times \R$ were invariant by some $1$-parameter group of rigid motions of the ambient space \cite{AbRo,MoOn,NeRo1,SaE,SaTo}. In connection with this theorem, we shall obtain an existence result for a \emph{Plateau problem at infinity} for complete surfaces with $H=1/2$ in $\H^2 \times \R$. We will also construct explicit examples of surfaces with $H=1/2$ and a prescribed simple hyperbolic Gauss map, and we will show that any surface with $H=1/2$ has a \emph{parallel-like} mean curvature one half surface with the same hyperbolic Gauss map. Also in Section 4 we will describe a Schwarz reflection principle in our context. In the end of the section, we shall see how the extended Lawson correspondence in \cite{Dan} lets us construct infinitely many complete minimal surfaces in the $3$-dimensional Heisenberg group with prescribed holomorphic quadratic differential.

At last, in Section 5 we turn our attention to minimal surfaces in $\H^2\times \R$. Given a complete minimal vertical graph in $\H^2\times \R$ with height function $h$, we can consider the \emph{canonical $1$-form} of the graph as the holomorphic part of $dh$. We will show that the map carrying each congruency class of complete minimal vertical graphs to its associated canonical $1$-form (defined with a $\pm$ ambiguity) is a bijection onto the space of holomorphic $1$-forms on $\C$ or $\D$, except for those of the form $\omega = c\, dz$ on $\C$ for some $c\in \C$. Our result also classifies all complete orientable minimal surfaces in $\H^2\times \R$ whose \emph{angle function} $u:\Sigma\flecha [-1,1]$ omits some interior value. In the end, we will apply some known results from the theory of harmonic maps into the Poincaré disk in order to describe domains in $\H^2$ over which a complete minimal vertical graph in $\H^2\times \R$ can be built or not.

\subsection*{Acknowledgements}

The authors wish to express their gratitude to the referee of this paper, whose remarks and questions have improved the exposition of the text, and have helped them to make the article more geometric.

Isabel Fernández was partially supported by MEC-FEDER Grant No. MTM2004-00160. Pablo Mira was partially supported by MEC-FEDER, Grant No. MTM2004-02746.

%%%%%

\section{Setup}\label{uno}

In this preliminary section we will describe some general facts that will be used in the study of both minimal surfaces and surfaces with $H=1/2$ in $\H^2\times \R$. First, we will analyze the structure and compatibility equations
of an immersed surface in $\H^2 \times \R$ in terms of a conformal
parameter for its first fundamental form. Subsequently, we will make some comments regarding harmonic maps into the hyperbolic plane $\H^2$ and their relation with spacelike CMC surfaces in the Lorentz-Minkowski $3$-space $\L^3$.

\subsection*{Integrability of surfaces in $\H^2\times \R$}

We shall realize $\H^2 \times \R$ in the Lorentz-Minkowski $4$-space
$\L^4$ as
$$\H^2 \times \R = \{(x_0,x_1,x_2,x_3)\in \L^4 : -x_0^2 +x_1^2 +x_2^2 =-1
, \hspace{0.3cm} x_0 >0\}.$$ Here we view $\L^4$ endowed with the
Lorentzian metric
$\esiz,\esde=-dx_0^2+dx_1^2+dx_2^2+dx_3^2,$ in canonical coordinates.

Let $\psi:\Sigma\flecha \H^2 \times \R$ be an immersed surface in $\H^2
\times \R$, that will be assumed to be simply connected. We will use the
notation $\psi =(N,h):\Sigma\flecha \H^2\times \R$ to denote the
\emph{vertical projection} $N:\Sigma\flecha \H^2$ and the \emph{height
function}
$h:\Sigma\flecha \R$ of $\psi$, respectively.

If $h$ is constant on an open set $\mathcal{U} \subset \Sigma$, then $\psi|_{\mathcal{U}}$ is a piece of a totally geodesic slice $\H^2 \times \{t_0\} \subset \H^2\times \R$. We shall rule out this trivial situation from now on, and consider only surfaces on which $h$ is never locally constant.

Let $\eta:\Sigma\flecha \S_1^3\subset \L^4$ denote the unit normal of
$\psi$ in $\H^2 \times \R$, where here $\S_1^3 =\{(x_0,x_1,x_2,x_3)\in
\L^4 : -x_0^2 +x_1^2 +x_2^2 +x_3^2=1 \} $ is the de Sitter $3$-space. In
this way, the metric conditions
$\esiz d\psi, \eta \esde = \esiz N ,\eta \esde =0$ hold, and the pair $\{\eta,N\}$ is
an orthonormal frame of the (Lorentzian) normal bundle of $\psi$ in
$\L^4$. We shall also use the splitting notation
$\eta = (\hat{N}, u):\Sigma\flecha \L^3 \times \R$. Following \cite{AlDa}, we shall call $u :\Sigma\flecha [-1,1]$ the \emph{angle function} of $\psi$.

As $\Sigma$ inherits via $\psi$ a Riemannian metric, it has an associated
Riemann surface structure. Thus, there exists a global conformal parameter
$z$ on $\Sigma$ with respect to which the induced metric of $\psi$ is
written as $\esiz d\psi, d\psi \esde = \landa |dz|^2 $ for a positive
smooth function $\landa$ on $\Sigma$. Let us also define the \emph{Hopf
differential} of $\psi$ as $p dz^2 = -\esiz \psi_z ,\eta_z\esde dz^2$,
i.e. as the $(2,0)$-part of its complexified second fundamental form. If
$H:\Sigma\flecha \R$ denotes the mean curvature of $\psi$ in $\H^2 \times
\R$ and $A:= -u h_z$, then it holds $$ \eta_z = -H \psi_z -
\frac{2p}{\landa} \psi_{\bar{z}} +A N.$$ More generally, if we consider
the moving frame
\begin{equation}\label{frame}
 \sigma=(\psi_z,\psi_{\bar{z}},\eta, N)^T
\end{equation}
the structure equations for the immersion are
 \begin{equation}\label{estruct}
  \sigma_z=\mathcal{U}
\sigma,\hspace{0.5cm} \sigma_{\bar{z}}=\mathcal{V}\sigma,
 \end{equation}
where \begin{equation*}
 \mathcal{U}=
\left(\def\arraystretch{1.3}\begin{array}{cccc}
 (\log \landa)_z & 0 & p & -h_z^2\\
 0 & 0& H \landa /2 & (\landa - 2|h_z|^2 )/2\\
 -H  & -2p/\landa & 0 & A \\
 1 - 2|h_z|^2/\landa & -2h_z^2/\landa & A & 0
\end{array}\right),
\end{equation*}
and
 \begin{equation*}
 \mathcal{V}= \left(\def\arraystretch{1.3}\begin{array}{cccc}
  0 & 0& H \landa /2 & (\landa - 2|h_z|^2 )/2\\
  0 & (\log \landa)_{\bar{z}} & \bar{p} & -h_{\bar{z}}^2 \\
 -2\bar{p}/\landa &  -H & 0 & \bar{A} \\
  -2h_{\bar{z}}^2/\landa & 1- 2|h_z|^2 /\landa & \bar{A} & 0
\end{array}\right).\end{equation*}
By examining the last coordinate in these equations, we deduce the
following relations between the coefficients:

\begin{equation}\label{lasces}
\left\{\def\arraystretch{1.3} \begin{array}{lccc} {\bf (C.1)} & h_{zz} & =
& \displaystyle \frac{\landa_z}{\landa} h_z + p u \\ {\bf (C.2)} &
h_{z\bar{z}} & = & \displaystyle \frac{\landa H}{2} u \\ {\bf (C.3)} &
u_{z} & = & - H h_z -\displaystyle \frac{2 p}{\landa} h_{\bar{z}}\\ {\bf
(C.4)} & \displaystyle \frac{4 |h_{z}|^2}{\landa} & = & 1 - u^2
\end{array}\right.
\end{equation}
The integrability condition of the system \eqref{estruct} is given by
\begin{equation}\label{entis}
\mathcal{U}_{\bar{z}} -\mathcal{V}_{z} +\left[\mathcal{U},\mathcal{V}\right] =0.
\end{equation}
Of this matrix identity, the entries $(1,1)$, $(1,3)$, $(1,4)$ and $(3,4)$ provide, respectively,
the Gauss-Codazzi-Ricci equations:

\begin{equation}\label{GCR}\def\arraystretch{1.6}
\begin{array}{lrcl} \text{Gauss: }& \landa(\log \landa) _{z\bar{z}}&=&
2\big(|p|^2 -\landa^2 (H^2 -1) /4 - \landa |h_z|^2 \big).\\
 \text{Codazzi: } & 2 p_{\bar{z}} &=& \landa ( H_z + A) \\
  \text{Codazzi (bis): } & -(h_z^2)_{\bar{z}} + (|h_z|^2 )_z & = & AH\landa /2 - \bar{A} p + \landa_z |h_z|^2 / \landa\\
 \text{Ricci: } & A_{\bar{z}} -\overline{A}_z &=& \displaystyle \frac{4i}{\landa}{\rm
 Im} (\overline{p} \, h_z^2).
\end{array}\end{equation}
All the other entries of \eqref{entis} provide relations that are trivial, or that follow from
\eqref{GCR}. Hence the equations \eqref{GCR} are the necessary and sufficient conditions for the
integrability of \eqref{estruct}.

Many of the equations we have obtained up to now are superfluous. For instance, the Ricci equation
follows directly from {\bf (C.3)}, and Codazzi (bis) is obtained by putting together {\bf (C.1)}
and {\bf (C.2)}. A somewhat lengthier computation also indicates that the Gauss equation may be
obtained from {\bf (C.1)}, ..., {\bf (C.4)}. At last, {\bf (C.1)} may be derived from {\bf (C.2)},
{\bf (C.3)} and {\bf (C.4)} by differentiation of {\bf (C.4)}. It is also important to observe
that, using {\bf (C.2)}, the Codazzi equation is written as
 \begin{equation}\label{cod}
 {\rm Codazzi:} \hspace{0.4cm} Q_{\bar z}=2pH_{\bar z} + \landa H H_z, \hspace{1cm}
Q:=2H p + h_z^2.
 \end{equation}
So, after these simplifications, only the Codazzi equation \eqref{cod} and the last three equations
in \eqref{lasces} remain as the integrability conditions of the system. More specifically, we have
proved:

\begin{pro}
Let $\Sigma$ denote a simply connected Riemann surface. The system \eqref{estruct} admits a
solution $\sigma :\Sigma\flecha \C^4\times \C^4 \times \L^4 \times \L^4$ if and only if the
coefficients $\landa,H,u,h :\Sigma\flecha \R$ and $p:\Sigma\flecha \C$ verify {\bf (C.2)}, {\bf
(C.3)}, {\bf (C.4)} and \eqref{cod}.
\end{pro}

This fact will be used in the proof of Theorem \ref{th:main} to show that by choosing an adequate
initial condition $\sigma (z_0)= \sigma_0$, the above equations actually produce a surface in
$\H^2\times \R$. In other words, equations {\bf (C.2)}, {\bf (C.3)}, {\bf (C.4)} and \eqref{cod}
are sufficient for the integrability of surfaces in $\H^2\times \R$. We do not detail this argument
here, as a proof of this last statement (actually of a more general one) is obtained in \cite{FeMi}
and in \cite{Dan}.

Moreover, we have the following consequence of \eqref{cod}:
 \begin{cor}[Abresch-Rosenberg]
The quadratic differential $Q dz^2$ is holomorphic on any surface with
constant mean curvature $H$ in $\H^2 \times \R$.
 \end{cor}

We shall refer to the holomorphic quadratic differential $Q dz^2$ as the
\emph{Abresch-Rosenberg differential} of a constant mean
curvature surface in $\H^2 \times \R$.

Let us also point out that the holomorphicity of $Q dz^2$ is not equivalent to the constancy  of
the mean curvature, as there surfaces whose Abresch-Rosenberg differential is holomorphic and that
are not CMC surfaces (see \cite{FeMi}).

\subsection*{Harmonic maps into the hyperbolic plane}

Let $\Sigma$ be an open simply connected Riemann surface, and let
$G:\Sigma\flecha \H^2$ be a smooth map. Then $G$ is harmonic if and only if the $(2,0)$-part of its complexified first fundamental form, i.e. $Q_0 dz^2 := \esiz G_z,G_z \esde
dz^2$, is a holomorphic quadratic differential. We shall call it the
\emph{Hopf differential} of the harmonic map $G$.

Let $\mu :\Sigma\flecha [0,+\8)$ be the smooth function so that
 \begin{equation}\label{afs}
\esiz dG,dG\esde = Q_0 dz^2 + \mu |dz|^2 + \bar{Q_0} d\bar{z}^2
 \end{equation}
holds for the harmonic map $G:\Sigma\flecha \H^2$. Then we have the following elementary facts:
 \begin{itemize}
 \item
As $\esiz dG,dG\esde$ is Riemannian, we get $\mu^2 - 4|Q_0|^2 \geq 0$, and equality holds exactly
at the singular points of $G$.
 \item
A point $z_0\in \Sigma$ is a \emph{branch point} of $G$ (i.e. $dG (z_0)=0$) if and only if $\mu (z_0)=Q_0(z_0)=0$.
 \item
$Q_0(z_0)=0$ at some point $z_0\in \Sigma$ if and only if $G$ is conformal (holomorphic or
antiholomorphic) at $z_0$. Moreover, any conformal map from $\Sigma$ into $\H^2$ is trivially a
harmonic map with $Q_0 =0$.
 \end{itemize}

Let us now relate the harmonic maps in $\H^2$ with surfaces in $\L^3$. For this, let us consider $f:\Sigma\flecha \L^3$ a spacelike surface in $\L^3$, oriented so that its unit normal $G$ takes its values in $\H^2$, i.e. in the upper sheet of the hyperboloid $\mathcal{H}^2 =\{x\in \L^3 : \esiz x,x\esde =-1\}.$ Let also $H:\Sigma\flecha \R$ denote its mean curvature. It is then well known that $G$ is a harmonic map into $\H^2$ if and only if $f$ is a CMC surface.

Suppose now that $H=1/2$, and write $\esiz df,df\esde= \tau_0|dz|^2$ for a positive smooth function
$\tau_0$. As $H=1/2$, the Hopf differential of $G$ agrees with the Hopf differential of $f$ in
$\L^3$. Furthermore, $\{Q_0,\tau_0\}$ verify the Gauss equation for $f$ in $\L^3$, that is,
 \begin{equation}\label{eq:gauss}
(\log \tau_0)_{z\bar{z}} = \tau_0 /8 -2|Q_0|^2 / \tau_0 .
 \end{equation}
In addition, the metric of the Gauss map $G$ is given by

 \begin{equation}\label{metG}
\esiz dG,dG\esde = Q_0 dz^2 + \mu |dz|^2 + \bar{Q_0} d\bar{z}^2, \hspace{0.5cm} \mu= \frac{\tau_0}{4} +
\frac{4|Q_0|^2}{\tau_0}.
 \end{equation}

Consider now the map $f^{\sharp} := f - 2G:\Sigma\flecha \L^3$, which is a parallel surface of $f$. It is well known, and also straightforward to check, that $f^{\sharp}:\Sigma\flecha \L^3$ has the same conformal structure that $f$, it has $G$ as its Gauss map, it has constant mean curvature $H^{\sharp} =-1/2$ and its Hopf differential is $Q^{\sharp} = -Q_0$.

Let $\tau^{\sharp}$ be the conformal factor of the metric of $f^{\sharp}$, i.e. $\esiz
df^{\sharp},df^{\sharp}\esde = \tau^{\sharp} |dz|^2$. The above comments imply that
$\{Q_0,\tau^{\sharp}\}$ verify the Gauss equation \eqref{eq:gauss}, and it also holds that
 \begin{equation}\label{mus}
\mu =\frac{\tau_0}{4} + \frac{4|Q_0|^2}{\tau_0} = \frac{\tau^{\sharp}}{4} + \frac{4
|Q_0|^2}{\tau^{\sharp}}.
 \end{equation}
Hence, $\tau^{\sharp} = 16 |Q_0|^2 /\tau_0$, what shows that the singular points of $f^{\sharp}$
are located at the umbilics of $f$.

Let us also observe that if we reverse the orientation of $f^{\sharp}$, we get a spacelike surface
in $\L^3$ with $H=1/2$ and whose Gauss map is $-G:\Sigma\flecha \H_{-}^2$, being $\H_{-}^2
=\mathcal{H}^2 \setminus \H^2$ the lower sheet of the hyperboloid $\mathcal{H}^2 \subset \L^3$.

Motivated by these facts, we formulate the following definition:

\begin{defi}\label{weidata}
Let $G:\Sigma\flecha \H^2$ be a harmonic map into $\H^2$ with Hopf differential $Q_0 dz^2$. We shall say that $G$ \emph{admits Weierstrass data} if there exists a smooth positive function $\tau_0:\Sigma\flecha (0,+\8)$ so that \eqref{metG} holds. In that case, the pair $\{Q_0,\tau_0\} $ will be called \emph{Weierstrass data} for $G$.
\end{defi}

The most obvious examples of harmonic maps admitting Weierstrass data are the Gauss maps of spacelike surfaces with $H=1/2$ in $\L^3$. This is just a consequence of the previous discussion. It is also immediate to realize that if $G:\Sigma\flecha \H^2$ is a conformal map without branch points, and we denote $\tau_0 := 8 \esiz G_z,G_{\bar{z}}\esde >0$, then $\{Q_0=0,\tau_0\}$ are Weierstrass data for $G$.

Conversely, we have the main conclusion of this subsection:
 \begin{lem}\label{lemarev}
If $\{Q_0,\tau_0\}$ are Weierstrass data for a harmonic map $G:\Sigma\flecha \H^2$, then they satisfy the Gauss equation \eqref{eq:gauss}.
 \end{lem}
\begin{proof}
If $G$ is conformal and regular, we have $Q_0 =0$, and \eqref{eq:gauss} holds trivially for $\tau_0 := 8 \esiz G_z,G_{\bar{z}}\esde >0$.

Now, assume that $G$ is not conformal. This means that $Q_0$ only has isolated zeros in $\Sigma$.
Let $\mathcal{Z}\subset \Sigma$ be the set of zeros of $Q_0$, and take $z\in \Sigma\setminus
\mathcal{Z}$. Then it is known (see \cite{AkNi}) that around $z$, $G$ is the Gauss map of a unique
(up to translations) spacelike surface $f$ in $\L^3$ with constant mean curvature $1/2$. Let $\tau$
denote the conformal factor of the metric of $f$. Then our previous discussion ensures that
$\{Q_0,\tau\}$ are Weierstrass data for $G$ around $z$, and that we have $\tau=\tau_0$ or $\tau =
16 |Q_0|^2 /\tau_0$. But at this point it is a direct computation to check that an arbitrary
positive smooth function $\delta$ on $\Sigma$ verifies \eqref{eq:gauss} with respect to $Q_0$ if
and only if $16 |Q_0|^2 /\delta$ verifies \eqref{eq:gauss}. Consequently, \eqref{eq:gauss} holds
for the original Weierstrass data $\{Q_0,\tau_0\}$ at every point $z\in \Sigma\setminus
\mathcal{Z}$. As $\mathcal{Z}$ is discrete, we conclude by continuity that $\{Q_0,\tau_0\}$ satisfy
\eqref{eq:gauss} globally on $\Sigma$. This concludes the proof.
\end{proof}

\noindent {\it Remark 1:} The above comments easily imply the following fact: if $G:\Sigma\flecha
\H^2$ is a harmonic map admitting Weierstrass data $\{Q_0,\tau_0\}$, and if $Q_0$ vanishes
somewhere, then $\tau_0$ is unique, i.e. the Weierstrass data are unique for $G$. In contrast, if
$Q_0$ never vanishes, then a second set of Weierstrass data, namely, $\{Q_0,
\tau^{\sharp}:=16|Q_0|^2 /\tau_0 \}$, is available for $G$. Moreover, by Lemma \ref{lemarev},
$\tau^{\sharp}$ still satisfies \eqref{eq:gauss}. Let us also observe that by \eqref{mus} we have
at every point in $\Sigma$ that
 \begin{equation}\label{yanise}
\{\tau_0,\tau^{\sharp}\} \in 2(\mu \pm \sqrt{\mu^2 - 4|Q_0|^2}).
 \end{equation}
So, $G$ is singular everywhere if and only if the solutions $\tau_0$ and $\tau^{\sharp}$ coincide.

\vspace{0.2cm}

We conclude this section with two elementary lemmas that will be used later on:

\begin{lem}\label{lem:Gsing}
Let $G:\Sigma\flecha \H^2$ be a harmonic map that is singular on an open subset of $\Sigma$. Then
$G(\Sigma)$ lies in a geodesic of $\H^2$.
\end{lem}
\begin{proof}
First of all, observe that by analyticity, $G$ is singular everywhere in $\Sigma$. As a result, we
have $\mu^2 = 4|Q_0|^2$ everywhere. It follows that if $Q_0$ vanishes identically, $G$ is constant
and the result trivially holds. If not, $Q_0$ only has isolated zeros. Take $z_0$ with $Q(z_0)\neq
0$. Then, by changing locally the complex parameter $z$ around $z_0$, we may assume that $Q_0 dz^2 = (1/4)
dz^2$, and so $\mu =1/2$. Denoting $z=s+it$, this implies by \eqref{metG} that $\esiz G_t,G_t\esde
= \esiz G_s,G_t\esde =0$, so $G=G(s)$. But finally, the harmonicity of $G$ indicates that $G_{ss}$
is collinear with $G$. In other words, $G$ parametrizes locally a piece of a geodesic in $\H^2$. By
analyticity, $G(\Sigma)$ lies in this geodesic and we are done.
\end{proof}

\begin{lem}\label{superre}
Let $G,\widetilde{G}:\Sigma\flecha \H^2$ be two harmonic maps for which \eqref{afs} holds for the
same functions $\mu,Q_0$. Choose $z_0\in \Sigma$ a regular point of $G$, and assume that
$G(z_0)=\widetilde{G}(z_0)$ and $G_z(z_0)=\widetilde{G}_z(z_0)$. Then $G=\widetilde{G}$.
\end{lem}
\begin{proof}
By the initial conditions, both $G,\widetilde{G}$ are regular surfaces in $\L^3$ around $z_0$.
Moreover, their respective orientations agree at $z_0$. By hypothesis, both surfaces have the same
first fundamental form. In addition, as both of them lie in $\H^2$ and their orientations agree,
their second fundamental forms also agree. Finally, as they share the initial conditions at $z_0$,
$G$ and $\widetilde{G}$ must coincide around $z_0$, and thus globally by analyticity.
\end{proof}

These two lemmas prove, in particular, that if $G:\Sigma\flecha \H^2$ is a harmonic map admitting
Weierstrass data $\{Q_0,\tau_0\}$, then any other harmonic map into $\H^2$ having the same
Weierstrass data $\{Q_0,\tau_0\}$ differs from $G$ just by an isometry of $\H^2$. In other words,
the Weierstrass data determine the harmonic map uniquely up to rigid motions.

%%%%%
\section{Surfaces of mean curvature one half in  $\H^2\times\R$}\label{tres}

Our aim in this section is to study the mean curvature one half surfaces in $\H^2\times \R$ in terms of an associated hyperbolic Gauss map. We will first of all describe this hyperbolic Gauss map, and then we will analyze how to recover a mean curvature one half surface from its hyperbolic Gauss map.

\subsection*{The hyperbolic Gauss map}

Let $\psi=(N,h):\Sigma\flecha \H^2 \times \R$ be an immersed surface whose
vertical projection $N:\Sigma\flecha \H^2$ is regular, i.e. $dN$ is a linear isomorphism at every point. By {\bf (C.4)}, the
regularity condition imposed to $N$ is equivalent to the fact that the angle function
$u:\Sigma\flecha [-1,1]$ never vanishes. This provides a canonical orientation for any surface of this type. Specifically, we will always assume that $\psi$ is oriented so that its angle function $u$ is positive.

So, for any surface in $\H^2\times \R$ with regular vertical projection
we can define the map
 \begin{equation*}
 \xi = \frac{1}{u} (\eta + N):\Sigma\flecha \N^3 :=\{x\in \L^4 : \esiz x,x\esde =0, x_0
 >0\}.
 \end{equation*}
Let us remark here that $\xi_0 >0$ because $N_0 >0$ and $\esiz \xi,N\esde =-1/u<0$.

If we observe that the last coordinate of $\xi$ is constantly $1$, the
fact that $\xi$ lies in $\N^3$ implies the existence of a map
$G:\Sigma\flecha \H^2\subset \L^3$ such that $\xi = (G,1):\Sigma\flecha
\L^3 \times \R \equiv \L^4$.

\begin{figure}[htbp]
\begin{center}
\includegraphics[width=.5\textwidth]{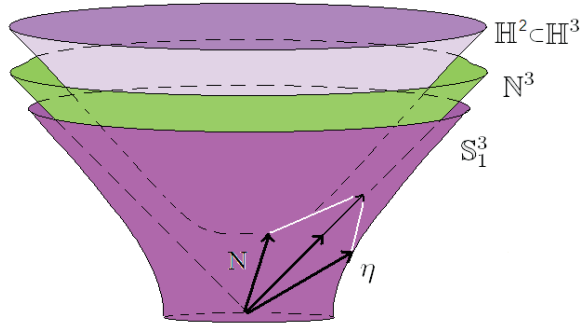}\caption{The normal vector $\xi$.}
\end{center}
\end{figure}

\begin{defi}
The map $G:\Sigma\flecha \H^2$ will be called the \emph{hyperbolic Gauss map} of the surface with regular vertical projection $\psi:\Sigma\flecha \H^2 \times \R$.
\end{defi}
The terminology of this definition has been motivated by its similarity with the construction of the classical hyperbolic Gauss map for surfaces in $\H^3$ \cite{Eps,Bry}. Specifically, let $f:\Sigma\flecha \H^3$ be a surface with unit normal $\nu:\Sigma\flecha \S_1^3$. Then the map $\cN:= f + \nu $ is a normal vector field taking its values
in the positive light cone $\N^3$. Moreover, if $\cN=(\cN_0, \dots, \cN_3)$, the map $\mathcal{G}:=(\cN_1/ \cN_0,\cN_2/ \cN_0,\cN_3/ \cN_0)$ can be viewed as a map from $\Sigma$ into $\S^2$, and satisfies that $\cN/\cN_0= (1,\mathcal{G})$. This map $\mathcal{G}$ is the hyperbolic Gauss map of $f$.

We also would like to remark that $G$ can indeed be called a Gauss map because of the following property: \emph{two surfaces in $\H^2\times \R$ with regular vertical projection meeting at a point $p\in \H^2\times \R$ have the same hyperbolic Gauss map $G(p)$ if and only they are tangent at $p$}. Indeed, two surfaces meeting at $p=(N,h)$ have the same hyperbolic Gauss map $G(p)$ if and only if their respective null normal vectors $\xi, \tilde{\xi}$ verify $\xi (p)=(G(p),1)= \tilde{\xi}(p)$, if and only if their respective normal bundles in $\L^4$ at $p$ are both spanned by $\{N,(G(p),1)\}$, if and only if they have the same tangent plane at $p$. We are grateful to the referee for this interesting observation.

The following theorem, whose proof follows easily from our previous
discussion, is the key tool for the main results of this work regarding surfaces with $H=1/2$ in $\H^2\times \R$.

\begin{teo}
The hyperbolic Gauss map of a mean curvature one half
surface in $\H^2\times \R$ with regular vertical projection is harmonic.
\end{teo}
\begin{proof}
A direct computation from the structure equations and {\bf (C.4)} shows that for every surface
$\psi =(N,h)$ in $\H^2\times \R$ it holds $$\esiz (\eta +N)_z,(\eta +N)_z \esde = (p+h_z^2 ) (2H-1
-u^2).$$ So, if $\psi$ has constant mean curvature $H=1/2$ and $Q$ is its Abresch-Rosenberg
differential, we see that
 \begin{equation}\label{temita}
\esiz \xi_z,\xi_z \esde = \frac{1}{u^2}\esiz (\eta +N)_z,(\eta +N)_z \esde = -Q.
 \end{equation}
As $Q$ is holomorphic, $\esiz \xi_{z\bar{z}} ,\xi_z \esde =0$, i.e. $\esiz G_{z\bar{z}} ,G_z\esde
=0$. So, $G$ is harmonic on $\H^2$.
\end{proof}

\vspace{0.2cm}

\noindent {\it Remark 2:} Let $\psi:\Sigma\flecha \H^2\times \R$ be a surface with mean curvature $H=1/2$ with respect to some selected orientation on $\Sigma$. Assume that $\psi$ has regular vertical projection. Then $H=\pm 1/2$ holds with respect to the canonical orientation on $\Sigma$ given by $u>0$. However, if $H=-1/2$ holds, its vertical symmetry is a surface in $\H^2\times \R$ with regular vertical projection for which $H=1/2$ holds for its canonical orientation. Thereby, we do not lose generality by restricting ourselves to the $H=1/2$ case.

\vspace{0.2cm}

It is immediate to observe from the above proof that if $\psi:\Sigma\flecha \H^2\times \R$ is a surface with $H=1/2$ and hyperbolic Gauss map $G$, then the Abresch-Rosenberg differential $Q$ of $\psi$ and the Hopf differential $Q_0$ of $G$ are related by $Q= - Q_0$. Moreover, we have

\begin{teo}\label{datW}
Let $G:\Sigma\flecha \H^2$ be the hyperbolic Gauss map of a mean curvature one half surface $\psi:\Sigma\flecha \H^2\times \R$ with regular vertical projection. Let $Q,\landa,u$ denote, respectively, the Abresch-Rosenberg differential, the conformal factor of the metric, and the angle function of $\psi$. Then $\{-Q,\landa u^2\}$ are Weierstrass data for $G$.

In particular, the hyperbolic Gauss map of a surface with $H=1/2$ always admits Weierstrass data.
\end{teo}
\begin{proof}
By the structure equations and \eqref{metG} we have $$\mu =2 \esiz G_z,G_{\bar{z}}\esde =\frac{2}{u^2} \esiz (\eta +N)_z,(\eta +N)_{\bar{z}} \esde = \frac{\landa u^2}{4} + \frac{4 |Q|^2}{\landa u^2}.$$ So, Definition \ref{weidata} and the relation $Q = -Q_0$ finish the proof.
\end{proof}
Let us also remark that by Theorem \ref{datW}, the hyperbolic Gauss map of a mean curvature one half surface with regular vertical projection cannot be constant.

\vspace{0.2cm}

\noindent {\it Remark 3:} A mean curvature one half surface in $\H^2 \times \R$ can have points with non-regular vertical projection, as evidenced by some helicoidal examples in \cite{SaTo}. But except for one particular example, these singular points have empty interior on a general mean curvature one half surface. Indeed, if $u=0$ holds on an open piece of a surface with $H=1/2$, a simple look at the integrability conditions shows that $Q$ vanishes, the height function $h$ is harmonic and the metric is flat. Thus we get a piece of a right cylinder over a horocycle in $\H^2$ (see \cite{AbRo,SaE}).

\vspace{0.2cm}

\noindent {\it Remark 4:} The harmonicity of this hyperbolic Gauss map has
close ties with the theory of Bryant surfaces in the hyperbolic
$3$-space $\H^3$. Let $f:\Sigma\flecha \H^3$ be a Bryant surface in $\H^3$ with unit normal $\nu :\Sigma\flecha \S_1^3$, and define again $\cN= f + \nu :\Sigma\flecha \N^3$. Then the $(2,0)$-part of the second fundamental form of $f$ in the direction of $\cN$ agrees with the
Hopf differential of the surface, $-\esiz f_z,\cN_z\esde =Q$, and thus is
holomorphic. In addition, the hyperbolic Gauss map $\mathcal{G}= \cN /
\cN_0$ turns out to be conformal \cite{Bry}.

In our situation of mean curvature one half surfaces in $\H^2 \times \R$,
the key observation is that the null normal vector field
$\eta +N$ also has the property that $-\esiz \psi_z, (\eta + N)_z \esde $
is holomorphic (it is the Abresch-Rosenberg differential). The final step is suggested by the case of Bryant surfaces, and
is to note that if we divide $\eta +N$ by its last coordinate (rather than by its first one), the map we obtain has its values in $\H^2$ and is harmonic.

\subsection*{Surfaces with prescribed hyperbolic Gauss map}

We have seen that any surface with $H=1/2$ and regular vertical
projection in $\H^2 \times \R$ has an associated harmonic Gauss map. Next
we shall deal with the inverse problem: \emph{can a harmonic map
$G:\Sigma\flecha \H^2$ from a simply connected open Riemann surface into $\H^2$ always be
realized as the hyperbolic Gauss map of a mean curvature one half surface in
$\H^2 \times \R$? }

The following lemma is a key ingredient in order to provide an answer to the above question.
 \begin{lem}\label{lem:A}
Let $\{-Q, 2\tau \} $ be Weierstrass data of a harmonic map from an
open simply connected Riemann surface $\Sigma$ into $\H^2$. Let $z_0\in
\Sigma$ be an arbitrary point in $\Sigma$. Then for any $\vt_0\in\C$ the second order differential
system
  \begin{equation}\label{eq:A}
\left\{ \begin{array}{lll}
h_{zz} & = & (\log\tau)_z\,h_z + Q\, \sqrt{\displaystyle\frac{\tau+2|h_z|^2}{\tau}}, \\
\\
h_{z\bar z} &=& \displaystyle\frac{1}{2}\sqrt{\tau(\tau+2|h_z|^2)},
\end{array}\right.
\end{equation}
has a  globally defined solution $h:\Sigma\flecha \R$, unique up to additive constants, satisfying the initial condition $h_z(z_0)=\vt_0$.
 \end{lem}

\begin{proof}
Consider firstly the first order differential system
 \begin{equation}\label{eq:Be}
\left\{ \begin{array}{lll}
\vt_z & = & (\log\tau)_z\,\vt + Q\, \sqrt{\displaystyle\frac{\tau+2|\vt|^2}{\tau}}, \\
\\
\vt_{\bar z} &=& \displaystyle\frac{1}{2}\sqrt{\tau(\tau+2|\vt|^2)}.
\end{array}\right.
\end{equation}
From (\ref{eq:gauss}) and $Q_{\bar z}=0,$ a straightforward computation indicates by means of \eqref{eq:A} that the following condition is satisfied:
\begin{displaymath}
\displaystyle\frac{\partial}{\partial \bar{z}} \left((\log\tau)_z\,\vt + Q \,\sqrt{\displaystyle\frac{\tau+2|\vt|^2}{\tau}}\right)
=
\displaystyle\frac{\partial}{\partial z} \left(
\displaystyle\frac{1}{2}\sqrt{\tau(\tau+2|\vt|^2)} \right).
\end{displaymath}
So, as $\Sigma$ is simply connected, by the Frobenius theorem we get the existence of a unique
global solution $\vt:\Sigma\flecha \C$ of \eqref{eq:Be} verifying the initial condition
$\vt(z_0)=\vt_0$. Observing now that by \eqref{eq:Be} it holds $\vt_{\bar{z}}\in\R$, we can finally
infer the existence of a function $h:\Sigma\flecha\R$, unique up to additive constants, satisfying
$h_z=\vt$.
\end{proof}

Now we can establish the main theorem of this section. For technical reasons, we will assume that the set of singular points of the harmonic map $G$ we start with has empty interior. The case in which $G$ is singular on an open set will be discussed in Section \ref{four}.

\begin{teo}\label{th:main}%%%aqui he quitado lo de que la G sea nowhere holmorphic...ademas he metido lo de que los puntos singulares tienen interior vacio
Let $G:\Sigma\flecha \H^2\subset \L^3$ be a harmonic map from an open simply connected Riemann surface into the hyperbolic plane admitting Weierstrass data $\{-Q,2 \tau\}$. Assume that the set of singular points of $G$ has empty interior, choose $z_0 \in \Sigma$ a regular point of $G$, and $\vt_0\in \C$.

Then, there is a unique (up to vertical translations) mean curvature one half surface $\psi =(N,h) :\Sigma\flecha \H^2 \times \R\subset \L^4$ such that:
 \begin{enumerate}
 \item[{\rm (i)}]
$G$ is the hyperbolic Gauss map of $\psi$.
 \item[{\rm (ii)}]
$\tau =\landa u^2 /2$, where $\esiz d\psi,d\psi\esde = \landa |dz|^2$ and $u$ is the angle function of $\psi$.
 \item[{\rm (iii)}]
$dh (z_0)= \vt_0 dz + \overline{\vt_0} d\bar{z}$.
 \end{enumerate}
Moreover, such $\psi$ can be recovered in terms of $G$ as
 \begin{equation}\label{rep}
\psi = \left( \displaystyle\frac{4\Re \Big(G_z \big(2\overline{Q}h_z + \tau h_{\bar{z}}\big)\Big)}{\tau^2-4|Q|^2}  + G\, \sqrt{\displaystyle\frac{\tau+2|h_z|^2}{\tau}},\, h\right).
 \end{equation}
Here $h:\Sigma\flecha \R$ is the unique (up to additive constants) solution to the differential
system \eqref{eq:A} with $h_z(z_0)=\vt_0$.
\end{teo}

Before coming to the proof of this result, let us make some brief comments regarding its hypothesis
and its most basic consequences:
 \begin{enumerate}
 \item
By Theorem \ref{datW},  the quantities $\{-Q,\landa u^2\}$ are always Weierstrass data for the hyperbolic Gauss map of a mean curvature one half surface in $\H^2\times \R$. So, if the Weierstrass data of $G$ in Theorem \ref{th:main} are unique (i.e. if $Q$ vanishes at some point), the condition {\rm (ii)} holds automatically on any surface with $H=1/2$ having $G$ as its hyperbolic Gauss map.
 \item
As explained in Section 2, the Weierstrass data of a harmonic map into $\H^2$ are defined up to the
ambiguity $\tau_0\leftrightarrow 16|Q_0|^2 /\tau_0$, if $Q_0$ never vanishes. If this is the case,
we obtain in Theorem \ref{th:main} two different surfaces with $H=1/2$ in $\H^2\times \R$ having
$G$ as hyperbolic Gauss map and with $dh (z_0)=\vt_0 dz + \overline{\vt_0} d\bar{z}$ (one for each
choice of $\tau_0$). A geometric interpretation of this duality will be carried out in Section
\ref{four}.
 \item
Equation \eqref{rep} is at first defined only at the regular points of $G$. However, as by
hypothesis the singular set of $G$ has empty interior, we obtain by continuity that \eqref{rep}
must actually hold everywhere in $\Sigma$.
 \item
It follows from Theorem \ref{datW}, Theorem \ref{th:main} and Proposition \ref{sae} (that will be
proved in Section \ref{four}) that a harmonic map from an open simply connected Riemann surface
$\Sigma$ into $\H^2$ is the hyperbolic Gauss map of a surface with $H=1/2$ in $\H^2\times \R$ if
and only if it admits Weierstrass data (see Definition \ref{weidata}).
 \item
It is directly deduced from Theorem \ref{th:main} that the space of mean curvature one half surfaces with the same hyperbolic Gauss map $G$ can be seen as a $3$-parameter family, where two parameters are obtained by varying the initial condition $\vt_0\in \C$, and the other one is the height $h(z_0)$. However, as the variation of this last parameter only produces a vertical translation on the surface, we will not take it into account. In other words, we will regard the space of surfaces with $H=1/2$ and the same hyperbolic Gauss map as a $2$-parameter family. We will discuss later on when two surfaces of this $2$-parameter family are actually congruent.
 \end{enumerate}

The remainder of this section is devoted to prove Theorem \ref{th:main}. So, we fix the notations of this theorem and consider the function $h:\Sigma\to\R$ given in Lemma \ref{lem:A} with $h_z(z_0)=\vt_0$. Now define  the following functions on $\Sigma$:

\begin{equation}\label{chanci}
u=\sqrt{ \displaystyle\frac{\tau}{\tau+2|h_z|^2} }\,,\qquad\quad
\lambda = 2\tau + 4|h_z|^2.
\end{equation}

\begin{lem}\label{lem:de1a4}
The functions $u,\landa,h$ verify the following identities on $\Sigma$:
\begin{eqnarray}
h_{zz}&=& (\log \lambda)_z\, h_z + (Q-h_z^2) u \label{eq:1},\\
h_{z\bar{z}}&=&\frac{\lambda u}{4} \label{eq:2},\\
u_z  &=& \frac{-1}{2}\,h_z-\frac{2Q-2h_z^2}{\lambda}\, h_{\bar{z}} =
\frac{-u^2}{2} h_z - \frac{2Q}{\lambda} h_{\bar z} \label{eq:3}, \\
u^2&=&1 - \frac{4 |h_z|^2}{\lambda}  \label{eq:4}.
\end{eqnarray}
\end{lem}
\begin{proof}
Formula \eqref{eq:4} is a direct consequence of the definition of the functions $u,h$ and $\landa$.
Moreover, \eqref{eq:1} and \eqref{eq:2} follow from \eqref{eq:A} and \eqref{chanci}.
At last, \eqref{eq:3} is obtained after differentiating $u^2$ in \eqref{chanci}, and with the aid
of \eqref{eq:A} and \eqref{eq:4}.
\end{proof}

This lemma shows that the functions $u,h$ and $\lambda$ satisfy the integrability equations ${\bf (C.1)}$ to ${\bf (C.4)}$ for $H=1/2$ and $p:=Q-h_z^2.$
In this way, by our general discussion in Section \ref{uno},  we infer that the system given by \eqref{estruct}
is integrable (recall that $Q$ is holomorphic). Thus, \eqref{estruct} has a globally defined solution $$\sigma=(\psi_z,\psi_{\bar z},\eta,N):\Sigma\to\C^4\times\C^4\times\L^4\times\L^4,$$
where $\psi:\Sigma\to\L^4.$ Moreover, $\psi$ is unique once we fix initial data $\sigma (z_0)=\sigma_0 \in \C^4\times\C^4\times\L^4\times\L^4$.

Our aim now is to check that for an adequate initial condition $\sigma_0$, the map $\psi$ is actually a regular surface in $\H^2\times\R,$ it has mean curvature one half, and its hyperbolic Gauss map coincides with $G.$ First, we define
\begin{equation}\label{sita}
\xi:=\frac{1}{u}(N+\eta).
\end{equation}
From \eqref{estruct}, \eqref{chanci}, \eqref{eq:3} and \eqref{eq:4} we obtain
\begin{equation}\def\arraystretch{2}\begin{array}{lll}\label{eq:xiz}
\xi_z  &= & \displaystyle\frac{u}{2}\, \psi_z - \frac{2Q}{\lambda u} \,\psi_{\bar z} + \frac{1}{2}\left(-u h_z + \frac{4 Q h_{\bar{z}}}{\landa u} \right)\xi \\
& =& \displaystyle\frac{u}{2}\, \psi_z - \frac{2Q}{\lambda u} \,\psi_{\bar z} - \frac{u}{2}\left(h_z -
\frac{2Q h_{\bar{z}}}{\tau} \right)\xi.
\end{array}\end{equation}
With this, we can fix initial conditions for $\sigma$ at $z_0\in \Sigma$ in the following way:

\begin{defi}[Initial data]\label{def:cond}
We will denote by $\sigma=(\psi_z,\psi_{\bar z},\eta,N)^T$ the unique solution of (\ref{estruct}) with the following initial conditions:
\begin{displaymath}
\begin{array}{l}
{\setlength\arraycolsep{0,5cm}
\begin{array}{ll}
\xi(z_0)=(G(z_0),1),\; & \xi_z(z_0)=(G_z(z_0),0),   \\[0,5cm]
N_3(z_0)=0,\; &   \langle  N,\xi\rangle (z_0)=-\displaystyle\frac{1}{u(z_0)}  ,
\end{array}} \\[1cm]
\langle  N,\xi_z\rangle (z_0)=\displaystyle\frac{1}{2}\left(h_z - \displaystyle\frac{2Q h_{\bar z}}{\tau} \right) (z_0) .
\end{array}\end{displaymath}
Here $z_0 \in \Sigma$ is a regular point of $G$ and $N_3$ denotes the last coordinate function of $N.$
\end{defi}

This definition guaranties that the vectors $\{\xi(z_0),\xi_z(z_0),\xi_{\bar z}(z_0)\}$ are linearly independent, as we have chosen $z_0$ as a regular point of $G$. Thus the value of $\psi_z(z_0)$ and $\psi_{\bar z}(z_0)$ can be computed in terms of the frame $\{\xi(z_0),\xi_z(z_0),\xi_{\bar z}(z_0)\}$ by means of \eqref{eq:xiz} and its conjugate expression.  Here, we have used that as $G$ is regular at $z_0$, by \eqref{metG} we have $\tau \neq 2|Q|$ at $z_0$. So, the conditions in Definition \ref{def:cond} determine completely the value of $\sigma(z_0)$, using \eqref{sita} and \eqref{eq:xiz}.

\begin{lem}\label{lem:cond0}
Let $\sigma :\Sigma\flecha \C^4\times\C^4\times\L^4\times\L^4$ be the map described in Definition \ref{def:cond}. Then the following metric relations hold on $z_0$:
\begin{displaymath}
\begin{array}{ll}
\langle  \psi_z,\psi_z\rangle |_{z_0}=0 ,& \langle  \psi_z,\psi_{\bar z}\rangle |_{z_0}={\lambda(z_0)}/{2}, \\[0,5cm]
\langle   N,N\rangle |_{z_0}=-1, & \langle  \eta,\eta\rangle |_{z_0}=1  ,\\[0,5cm]
\langle  N,\psi_z\rangle |_{z_0}=0,  & \langle  \eta,\psi_z\rangle |_{z_0}=0 ,\\[0,5cm]  \langle  N,\eta\rangle |_{z_0}=0   .
\end{array}
\end{displaymath}
\end{lem}
\begin{proof}
We omit the point $z_0$ for the sake of simpleness.
Start by noting that, by \eqref{eq:xiz}, $$  0= \langle  \xi_z,\xi\rangle  = \frac{u}{2}\langle  \psi_z,\xi\rangle  -\frac{2Q}{u\lambda} \langle  \psi_{\bar z},\xi\rangle,  $$
and so $\langle  \psi_z,\xi\rangle =0.$ In addition, the following relations follow from \eqref{eq:xiz}, \eqref{metG} and Definition \ref{def:cond}:

{\setlength\arraycolsep{2pt}
\begin{eqnarray*}
-Q = \esiz G_z,G_z\esde &= \langle  \xi_z,\xi_z\rangle  =& \frac{u}{2}\langle  \psi_z,\xi_z\rangle  -\frac{2Q}{\lambda u} \langle  \psi_{\bar z},\xi_z\rangle  .
\\[0,5cm]
 \frac{1}{4}\left(\tau+\frac{4|Q|^2}{\tau}\right)=  \esiz G_z,G_{\bar{z}} \esde & = \langle  \xi_{\bar z},\xi_{z}\rangle = &
\frac{u}{2} \langle  \psi_{\bar z},\xi_{z}\rangle  -\frac{2\overline{Q}}{\lambda u} \langle  \psi_{z},\xi_{z}\rangle .
\end{eqnarray*}}
Hence, ${\langle  \psi_z,\xi_z\rangle = -Q/u}$ and  ${\langle  \psi_{\bar z},\xi_z\rangle =\lambda u/4}.$ This lets us write
{\setlength\arraycolsep{2pt}
\begin{eqnarray*}
 -\frac{Q}{u} &=\langle  \psi_z,\xi_z\rangle  = &
\frac{u}{2} \langle  \psi_z,\psi_z\rangle   -\frac{2Q}{\lambda u}\langle  \psi_z,\psi_{\bar z}\rangle, \\[0,5cm]
 \frac{\lambda u}{4} &=\langle  \psi_{z},\xi_{\bar{z}}\rangle  =& \frac{u}{2} \langle  \psi_{z},\psi_{\bar{z}}\rangle  -
\frac{2\overline{Q}}{\lambda u} \langle  \psi_{z},\psi_{z}\rangle ,
\end{eqnarray*}}
what gives
${ \langle  \psi_z,\psi_z\rangle }=0$ and
${\langle  \psi_{ z},\psi_{\bar z}\rangle =\lambda/2}.$

Let us see now that $\langle  N,N\rangle |_{z_0}=-1.$ Indeed, write $N(z_0)=(n,0)\in\L^3\times\R=\L^4.$
Since $z_0$ is a regular point for $G$ we can write in $\L^3$
 \begin{equation}\label{paran}
n=\alpha\,G_z+\bar\alpha \,G_{\bar z}+\gamma G, \qquad\alpha,\gamma\in\R .
 \end{equation}
These numbers $\alpha$ and $\gamma$ can be computed from the conditions in Definition \ref{def:cond}, to obtain

 \begin{equation}\label{paran2}
\alpha=2\,  \frac
{{2\,\overline{Q}h_z} +\tau h_{\bar{z}}}
{\tau^2-4|Q|^2}
\qquad \mbox{and}\qquad
\gamma=1/u .
 \end{equation}
Therefore
\begin{equation*}
\langle  N,N\rangle =\langle  n,n\rangle =-\alpha^2 Q- \bar\alpha^2 \overline{Q}  + \frac{2|\alpha|^2}{4}\left(\tau+\frac{4|Q|^2}{\tau}\right)
-\frac{1}{u^2}=-1.
\end{equation*}
From this and Definition \ref{def:cond} we also get
$${\langle  N,\eta\rangle |_{z_0}=0,\quad \langle  \eta,\eta\rangle |_{z_0}=1}.$$
Furthermore,
 \begin{equation}\label{hisa}
 \frac{h_z}{2} - \frac{Q h_{\bar{z}}}{\tau} =\esiz G_z,N\esde= \langle
\xi_z,N\rangle = \frac{u}{2}\langle  \psi_z,N\rangle -\frac{2Q}{u\lambda} \langle  \psi_{\bar
z},N\rangle  +\frac{h_z}{2} - \frac{Q h_{\bar{z}}}{\tau},
 \end{equation}
from where ${\langle \psi_z,N\rangle =0}.$ Finally, by $\esiz \psi_z,\xi\esde =0$, this leads to
${\langle \psi_z,\eta\rangle =0}.$
\end{proof}

\begin{proof1}
We will start proving existence, by showing that the map $\psi:\Sigma\flecha \L^4$ obtained via \eqref{estruct}, \eqref{chanci} and Definition \ref{def:cond} is a regular surface in $\H^2 \times \R$ with the desired conditions.

First of all, we use a standard argument to ensure that the metric relations appearing in Lemma \ref{lem:cond0} actually hold everywhere on $\Sigma$, and not just at $z_0$.

Let us introduce the notation $\sigma=(\psi_z,\psi_{\bar z},\eta,N)^T=(\sigma_1,\sigma_2,\sigma_3,\sigma_4)^T.$
From \eqref{estruct} we can deduce that the functions $\Phi_{i,j}:=\langle  \sigma_i,\sigma_j\rangle ,$ $i,j=1,\ldots,4,$
satisfy the following linear system of partial differential equations,
$$\left\{ \begin{array}{ccl}
(\Phi_{i,j})_z&=& \sum_{k=1}^4 \Big( U_{i,k}\Phi_{k,j} + U_{j,k} \Phi_{k,i}\Big),\\
(\Phi_{i,j})_{\bar z}&=& \sum_{k=1}^4 \Big( V_{i,k}\Phi_{k,j} + V_{j,k} \Phi_{k,i}\Big),
\end{array}\right.$$
where $U_{i,j}$ (resp. $V_{i,j}$) denotes the $(i,j)$ element of the matrix $\mathcal U$ (resp. $\mathcal V$).

On the other hand, it is direct to see that the functions $\phi_{i,j}=\phi_{j,i}$ described by
$$\phi_{1,1}=\phi_{2,2}=\phi_{1,3}=\phi_{2,3}=\phi_{1,4}=\phi_{2,4}=\phi_{3,4}=0,$$
$$\phi_{1,2}=\lambda/2, \qquad \phi_{3,3}=-\phi_{4,4}=1$$
also provide a solution to this system.
Since both solutions coincide at $z_0$ (Lemma \ref{lem:cond0}) they must coincide at any point of $\Sigma.$ Particularly, $\psi$ is a regular spacelike immersion into $\L^4$.

Now let us show that
 \begin{equation}\label{cambio1}
(\psi_z)_3=h_z,\hspace{0.5cm}  N_3=0 \hspace{0.5cm} \text{ and } \hspace{0.5cm} \eta_3=u.
 \end{equation}
Let $\gamma:=((\psi_z)_3,(\psi_{\bar z})_3,\eta_3,N_3)^T.$ Then
$$\gamma_z=\mathcal{U}\gamma\quad \mbox{ and }\quad \gamma_{\bar{z}}={\mathcal V}\gamma.$$
This system also has $(h_z,h_{\bar z},u,0)$ as a solution, by {\bf (C.1)} to ${\bf (C.4)}$. So, we only need to check that both solutions agree at $z_0.$ By the conditions
$N_3(z_0)=0$ and $\xi_3(z_0)=1$ of Definition \ref{def:cond}, we obtain $\eta_3(z_0)=u(z_0).$
In addition, since $(\xi_z)_3(z_0)=0,$ by examining the last coordinate of \eqref{eq:xiz} we get $(\psi_z)_3(z_0)=h_z(z_0),$ as we wished to show.

Let us check next that $\psi$ lies in $\H^2 \times \R$.
By \eqref{estruct} we know that
$$N_z=\left(1-\frac{2|h_z|^2}{\lambda}\right)\psi_z-\frac{2h_z^2}{\lambda}\psi_{\bar z} -u h_z\, \eta.$$
On the other hand, \eqref{cambio1} gives
$$\tilde{h}_z=\frac{2|h_z|^2}{\lambda}\psi_z+\frac{2h_z^2}{\lambda}\psi_{\bar z}+u h_z \,\eta ,$$ where here we are denoting $\tilde{h}=(0,0,0,h).$ More specifically, this identity follows from \eqref{cambio1}, since the equations in \eqref{cambio1} imply the relations $\esiz \tilde{h}_z,\eta \esde = u h_z$, $\esiz \tilde{h}_z,N\esde =0$, $\esiz \tilde{h}_z,\psi_z \esde=h_z^2$ and $\esiz \tilde{h}_z,\psi_{\bar{z}}\esde = |h_z|^2$.

With this, $\psi_z=N_z+\tilde{h}_z=(N_z,h_z)$ and therefore, up to a translation, $\psi=(N,h)$. So, $\psi$ lies in $\H^2\times \R$, and by \eqref{estruct} it has mean curvature one half, and its unit normal is $\eta$. Moreover, by its own construction, it holds $\tau = \landa u^2 /2$ and $dh(z_0)=\vt_0 dz + \overline{\vt_0}d\bar{z}$. It is also straightforward that $\psi$ has regular vertical projection.

It remains to check that the hyperbolic Gauss map of $\psi$ is $G.$
Let $\cG$ denote the hyperbolic Gauss map of $\psi$, that is, $\xi =(\mathcal{G},1)$.
From \eqref{eq:xiz} we have
$$\langle  \cG_z,\cG_z\rangle =-Q=\langle  G_z,G_z\rangle $$ and
$$ \langle  \cG_{z},\cG_{\bar z}\rangle =\frac{1}{4}\left(\tau+\frac{4|Q|^2}{\tau}\right)=\langle  G_z,G_{\bar z}\rangle. $$ As by Definition
 \ref{def:cond} we have
$\cG(z_0)=G(z_0)$ and $\cG_z (z_0)=G_z (z_0)$, and as $z_0$ is a regular point of $G$, we conclude as in Lemma \ref{superre} that $\cG= G$ on $\Sigma$.
This finishes the existence part.

For uniqueness, suppose that $\tilde{\psi}:\Sigma\to\H^2\times\R,$
$\tilde{\psi}=(\tilde{N},\tilde{h}),$ is another mean curvature one half surface in the conditions
of the theorem, and let $\tilde{\landa}$ denote the conformal factor of its metric. Let
$\tilde{\eta}$ denote its unit normal, with the last coordinate $\tilde{u}>0$, and assume that
$dh(z_0)=d\tilde{h}(z_0)$. Then by $\tau=\tilde{\lambda}\tilde{u}^2/2$ and {\bf (C.1)} to {\bf
(C.4)} for $\tilde{\psi}$, a somewhat lengthy but otherwise direct computation  shows that
$\tilde{h}_z$ satisfies the differential system \eqref{eq:A}. Consequently $\tilde{h}_z=h_z$.

Moreover, by hypothesis we have $\lambda u^2/2=\tau=\tilde{\lambda}\tilde{u}^2/2.$ But using {\bf (C.4)} for $\tilde{\psi}$ it is easy to check that
$$ \tilde{u}=\sqrt{ \displaystyle\frac{\tau}{\tau+2|\tilde{h}_z|^2} }\,= \sqrt{ \displaystyle\frac{\tau}{\tau+2|{h}_z|^2} }\, =u ,$$
and so it also holds $\tilde{\lambda}=\lambda$.

Thus the structure equations for $\tilde{\psi}$ and $\psi$ are the same, and their respective moving frames coincide at $z_0$. This implies that $\tilde{\psi}=\psi$ up to vertical translations in $\H^2 \times \R$.

Finally, the expression \eqref{rep} follows directly by repeating the computations described in \eqref{paran} and \eqref{paran2} at an arbitrary regular point $z\in \Sigma$, and not just at $z_0$, and using that the set of regular points of $G$ is dense by hypothesis. This ends up the proof.
\end{proof1}

\vspace{0.2cm}

Let $\psi=(N,h):\Sigma\flecha\H^2\times\R$ be a surface with constant mean curvature one half and
regular vertical projection, and suppose that its hyperbolic Gauss map $G$ is in the conditions of
Theorem \ref{th:main}. As a consequence of this theorem, up to vertical translations, $\psi$ is
uniquely determined (once we fix Weierstrass data for $G$) by the value of the differential of its
height function $h$  at a fixed point $z_0\in\Sigma$. Moreover, from \eqref{hisa} we infer that, at
the regular points of $G$, we have
\begin{equation}\label{enfins}
h_z=\displaystyle\frac{2 \tau\, ( \tau\langle N,G_z\rangle +2Q\langle N,G_{\bar z}
\rangle)}{\tau^2-4|Q|^2}.
\end{equation}
This shows that $\psi=(N,h)$ is uniquely determined (up to vertical translations, and after prescribing
$G$ and Weierstrass data for $G$) also by the value of its vertical projection $N(z_0)$ at an
arbitrary regular point $z_0\in \Sigma$ of $G$.

Conversely, the quantity $\psi (z_0)$ can be specified as initial condition for recovering a mean
curvature one half surface in terms of its hyperbolic Gauss map. Indeed, we have

\begin{cor}\label{co:main}
Let $G:\Sigma\flecha \H^2$ be a harmonic map from an open simply connected Riemann surface into the
hyperbolic plane admitting Weierstrass data $\{-Q,2 \tau\}$. Assume that the set of singular points
of $G$ has empty interior, choose $z_0 \in \Sigma$ a regular point of $G$, and $\psi_0\in
\H^2\times\R$.

Then, there is a unique mean curvature one half surface $\psi =(N,h) :\Sigma\flecha \H^2 \times \R$ such that:
 \begin{enumerate}
 \item[{\rm (i)}]
$G$ is the hyperbolic Gauss map of $\psi$.
 \item[{\rm (ii)}]
$\tau =\landa u^2 /2$, where $\esiz d\psi,d\psi\esde = \landa |dz|^2$ and $u$ is the angle function of $\psi$.
 \item[{\rm (iii)}]
$\psi(z_0)= \psi_0$.
 \end{enumerate}
\end{cor}

\begin{proof} Let $\psi_0=(N_0,h_0)\in\H^2\times\R$ and
consider, motivated by \eqref{enfins}, the complex number $\vt_0\in\C$ given by
\begin{equation}\label{eq:psi0}
 \vt_0=\displaystyle\frac{\,2\alfa_0 \tau^2(z_0)+ 4\tau (z_0)\overline{\alfa_0}Q(z_0)\,}{\tau^2(z_0)-4|Q(z_0)|^2},
\end{equation}
where $\alfa_0=\langle N_0,G_z(z_0)\rangle\in\C$ (recall that $z_0$ is a regular point of $G$ and
therefore $\tau^2-4|Q|^2\neq 0$ at $z_0$). It is obvious that from \eqref{eq:psi0} we actually have
that
 \begin{equation}\label{masfins}
\alfa_0 =\frac{\vt_0}{2} -\frac{Q \bar{\vt_0}}{\tau}.
 \end{equation}
Let $\psi=(N,h):\Sigma\to\H^2\times\R$ be the mean curvature one half surface obtained from Theorem
\ref{th:main} in terms of  $G$, $\{-Q,2\tau\}$ and $\vt_0$, and with $h(z_0)=h_0$. By \eqref{rep}
and \eqref{masfins} it is immediate to check that $$\esiz N(z_0),G_z (z_0)\esde=\frac{\vt_0}{2}
-\frac{Q \bar{\vt_0}}{\tau} = \alfa_0 =\langle N_0,G_z(z_0)\rangle .$$ But now, as $z_0$ is a
regular point of $G$, and we also know that $N_0\in \H^2$, $N(z_0)\in \H^2$, we must necessarily
have $N(z_0)= N_0$. This completes the existence part. Uniqueness was already proved above.
\end{proof}

%%%%%%%%%%%%%%%%%%%%%%%%%%%%%%%%%%%%%%%%%%%%%%%%%%%%%%%%%%%%%%%%%%%%%%%%%%%%%%%%%%
\vspace{0.2cm}

\noindent {\it Remark 5:} The above corollary has the following consequence: given an open simply connected Riemann surface $\Sigma$, the space of conformal immersions of $\Sigma$ into $\H^2\times \R$ with $H=1/2$ and regular vertical projection can be parametrized in terms of the space of harmonic maps from $\Sigma$ into the Poincaré disk that admit Weierstrass data, up to a certain initial condition.

\vspace{0.2cm}

Let $G:\Sigma\flecha \H^2$ be a harmonic map in the conditions of Corollary \ref{co:main}, and let
$\{-Q,2\tau\}$ be Weierstrass data for $G$. Then, by Corollary \ref{co:main}, it follows that the
class of simply connected mean curvature one half surfaces with $G$ as hyperbolic Gauss map is a
two-parameter family, where the parameters are given by the variation of the initial condition
$\psi_0\in \H^2\times\R$ (more specifically, of $N_0\in \H^2)$. Then, a natural question arising
here is whether all these surfaces are mutually non-congruent or not. This will depend on the
symmetries of $G$, as we discuss next.

Let $\psi,\tilde{\psi}:\Sigma\flecha \H^2\times \R$ be two surfaces with $H=1/2$ and the same hyperbolic Gauss map $G:\Sigma\flecha\H^2$, and assume that they are congruent. Thus there is some positive rigid motion $\Phi$ of $\H^2\times \R$ and some automorphism $\Gamma$ of the Riemann surface $\Sigma$ such that $\Phi \circ \psi = \tilde{\psi} \circ \Gamma$.  As both $\psi$ and $\tilde{\psi}$ are canonically oriented as surfaces with regular vertical projection (i.e. their angle functions are positive), we find that the positive rigid motion $\Phi$ must preserve the orientation of the vertical factor. So, $\Phi = (\Psi, {\rm Id} +c)$, where $c\in \R$ and $\Psi:\H^2 \subset \L^3\flecha \H^2\subset \L^3$ is an isometry of $\H^2$. Then we get that $\Psi\circ G = G \circ \Gamma$. That is, $\Psi$ is a symmetry of the hyperbolic Gauss map $G$.

Conversely, let $G:\Sigma\flecha \H^2$ denote a regular harmonic map on an open simply connected
Riemann surface $\Sigma$ admitting Weierstrass data $\{-Q,2\tau\}$, and let $\Psi:\H^2\flecha \H^2$
be a symmetry of $G$. That is to say, $\Psi \circ G =G\circ \Gamma$ for some automorphism  $\Gamma$
of $\Sigma$. Let $\psi:\Sigma\flecha \H^2\times \R$ be a mean curvature one half surface with
hyperbolic Gauss map $G$ constructed as in Theorem \ref{th:main} (or as in Corollary
\ref{co:main}). It is then clear that $\widetilde{\psi} := \Psi^{-1}\circ \psi \circ \Gamma
:\Sigma\flecha \H^2\times \R$ is the surface with $H=1/2$ and hyperbolic Gauss map $G$ constructed
via Corollary \ref{co:main} from $Q,\tau$ and $\widetilde{\psi}_0 := \Psi^{-1} (\psi (\Gamma
(z_0))) \in \H^2\times \R$.

All of this draws the following conclusions:

\begin{itemize}
\item
Given a harmonic map $G:\Sigma\flecha \H^2$ admitting Weierstrass data, if two initial conditions $\psi_0=(N_0,h_0),\tilde{\psi_0}=(\tilde{N_0},h_0)\in \H^2\times\R$ determine via Corollary \ref{co:main} congruent mean curvature one half surfaces, then $G$ necessarily has some symmetry. Thus, in the generic case, the class of non-congruent simply connected surfaces with $H=1/2$ sharing the hyperbolic Gauss map is a continuous two-parameter family.
 \item
If we fix $\psi_0=(N_0,h_0)\in \H^2\times\R$ an arbitrary point, each symmetry of $G$  provides a
point $\tilde{\psi}_0=(\tilde{N_0},h_0)\in \H^2\times\R$ such that $\{G,\tilde{\psi_0}\}$ and
$\{G,\psi_0\}$ generate congruent surfaces with $H=1/2$. In this sense, the larger is the symmetry
group of $G$, the smaller is the class of non-congruent mean curvature one half surfaces with $G$
as hyperbolic Gauss map.
\end{itemize}

\vspace{0.2cm}

\noindent {\it Remark 6:} Theorem \ref{th:main} indicates clearly that the theory of surfaces with
$H=1/2$ in $\H^2 \times \R$ is an appropriate setting for applying integrable systems techniques.
Indeed, one can associate to a harmonic map into $\H^2$ the whole machinery of a spectral
parameter, a zero-curvature representation, a Sym-Bobenko formula, and several Backlund-Darboux
transformations. These transformations let us construct via Theorem \ref{th:main} new mean
curvature one half surfaces in $\H^2\times \R$ starting from a previously known example. Let us
also observe that away from the zeros of the Abresch-Rosenberg differential, and up to a conformal
reparametrization, by \eqref{eq:gauss} the local geometry of a surface with $H=1/2$ in $\H^2\times
\R$ is modelled by the elliptic Sinh-Gordon equation: $\Delta \tau_0 = \sinh (\tau_0)$, where
$\Delta$ is the Euclidean Laplacian.

\section{Applications}\label{four}

We will describe in this section how the hyperbolic Gauss map we have introduced can be used to investigate global properties of mean curvature one half surfaces in $\H^2 \times \R$, and especially, to construct complete examples of such surfaces.

\subsection*{Existence of complete examples}

As a surprisingly simple consequence of our discussion, we can conclude that there are no restrictions in prescribing the Abresch-Rosenberg differential for \emph{complete} mean curvature one half surfaces in $\H^2 \times \R$:

\begin{teo}\label{exicom}
Any holomorphic quadratic differential on an open simply connected Riemann surface $\Sigma$ is the Abresch-Rosenberg differential of some complete surface with $H=1/2$ in $\H^2\times \R$. Moreover, the space of non-congruent complete mean curvature one half surfaces in $\H^2\times \R$ with the same Abresch-Rosenberg differential is generically infinite.
\end{teo}
\begin{proof}
Observe first of all that the right vertical cylinder over a horocycle in $\H^2\times \R$ has parabolic conformal structure (as it is flat and complete), and vanishing Hopf differential. So, we only need to consider the case in which $\Sigma\equiv \D$, or $\Sigma\equiv \C$ and $Q\not\equiv 0$.

Given a holomorphic quadratic differential $Q dz^2$ as above, Wan and Au \cite{Wan,WaAu} showed that there exists a unique surface $f:\Sigma\flecha \L^3$ with $H=1/2$ whose Hopf differential is $-Q dz^2$, and whose induced metric $\esiz df,df\esde = \tau_0 |dz|^2$ is complete. Let $\nu:\Sigma\flecha \H^2 \cup \H_{-}^2 \subset \L^3$ denote the Gauss map of $f$, and consider the rigid motion $P$ of $\L^3$ given by $P={\rm Id}$ if $\nu (\Sigma)\subset \H^2$ and $P(x_0,x_1,x_2)=(-x_0,x_1,x_2)$ if $\nu (\Sigma)\subset \H^2_{-}$. Here $\H^2_{-} = \{(-x_0,x_1,x_2) \in \L^3 : (x_0,x_1,x_2)\in \H^2\subset \L^3\}$. Let $G:=P\circ \nu:\Sigma\flecha \H^2$, and let $\psi:\Sigma\flecha \H^2 \times \R$ denote any of the mean curvature one half surfaces constructed from $G$ and $\tau_0$ via Theorem \ref{th:main}. Recall here that in the generic case, the family of such mean curvature one half surfaces is $2$-parametric by our discussion at the end of Section \ref{tres}. By its construction, the Abresch-Rosenberg differential of $\psi$ is precisely $Q dz^2$. In addition, we have $\esiz d\psi,d\psi\esde =\landa |dz|^2$, where
 \begin{equation}\label{desmetr}
 \tau_0 = \landa u^2 \leq \lambda .
  \end{equation}
Thus, by the completeness of $f$ we can conclude the completeness of $\psi$.
\end{proof}

Following the path suggested by this theorem, let us formulate the following \emph{Plateau problem at infinity}: Let $\gamma:\S^1 \flecha \S^1$ be a continuous homeomorphism. \emph{Is there a complete mean curvature one half surface $\psi:\D \flecha \H^2 \times \R$ whose hyperbolic Gauss map $G:\D\flecha \D$ extends continuously to $\overline{\D}$ and verifies $G|_{\S^1} = \gamma$ ?} Here $\D$ is the unit disk and $\H^2$ has been identified with the Poincaré disk $(\D,ds^2)$.

The following result follows immediately from \eqref{desmetr} and \cite{Aku,LiTa1,LiTa2}, where the Dirichlet problem at infinity for harmonic maps $G:\D \flecha \H^2$ is solved.

\begin{teo}
If $\gamma:\S^1 \flecha \S^1$ is a $C^{1,\alfa}$-diffeomorphism, $0<\alfa<1$, with ${\rm deg} (\gamma)=1$, then the above Plateau problem at infinity for mean curvature one half surfaces has at least a solution, and generically an infinite number of them.
\end{teo}

Apart from the previous ones, there are many results on the global construction of harmonic maps into the Poincaré disk. By our discussion, all these existence results translate directly into global existence results for mean curvature one half surfaces in $\H^2 \times \R$ by means of Theorem \ref{th:main}. It seems an interesting problem to analyze if the constructions of this paper and the results from the theory of harmonic maps can be applied to solve the Bernstein problem for mean curvature one half surfaces, i.e. to find all the entire vertical graphs with $H=1/2$ in $\H^2 \times \R$.

Let us also mention the following important open problem in the theory of harmonic maps, due to Schoen \cite{Sch,ScYa}: \emph{are there global harmonic diffeomorphisms from the complex plane onto the Poincaré disk?} This problem has been widely investigated by means of the related theory of spacelike CMC surfaces in $\L^3$. In this sense, building a geometric theory with an associated harmonic Gauss map into the Poincaré disk is interesting, since it may help to achieve a solution to the above problem. Indeed, in our situation, this question can be formulated as follows: \emph{is there a mean curvature one half surface in $\H^2\times \R$ with regular vertical projection and parabolic conformal structure, and whose hyperbolic Gauss map is a global diffeomorphism?}

\subsection*{Surfaces with singular hyperbolic Gauss map}

In Theorem \ref{th:main} we avoided the consideration of the case in which the harmonic map $G$ is singular on an open set. The next result deals with this remaining case.

\begin{pro}\label{sae}
Let $\psi:\Sigma\flecha \H^2\times \R$ be a mean curvature one half surface whose hyperbolic Gauss map $G:\Sigma\flecha \H^2$ is singular on an open set of $\Sigma$. Then $G(\Sigma)$ lies on a geodesic of $\H^2$, and $\psi$ is one of Sa Earp's standard hyperbolic screw motion examples.
\end{pro}
\begin{proof}
Let $\{-Q,2\tau\}$ be Weierstrass data for $G$. By Lemma \ref{lem:Gsing}, $G$ parametrizes a piece of a geodesic in $\H^2$.

Let us find the surface $\psi$ explicitly in these conditions. For that, we assume that $$G=G(t)=
(\cosh(t),\sinh(t),0):\R\flecha \H^2\subset \L^3,$$  where here $z=s+it$ is a global conformal
parameter for $\psi$. In particular, $2Q =\tau = 1/2$, and thus $\lambda u^2=2\tau=1$. In these
conditions, the system \eqref{eq:A} turns into
$$h_{zz} = h_{z\bar{z}}= \frac{1}{4}\sqrt{1+4|h_z|^2},$$
or equivalently
 \begin{equation*}%\label{siscutre}
 \left\{\def\arraystretch{1.4} \begin{array}{lll}
h_{ss} &=& \sqrt{1+h_s^2+h_t^2} \\
h_{st} &=& 0\\
h_{tt} &=& 0.
 \end{array} \right.
 \end{equation*}
This system can be explicitly integrated, to obtain
 \begin{equation}\label{laache}
h(s,t)=\left(\sqrt{1+y^2}\right)\, \cosh(s+s_0)+y\,t + c,
 \end{equation}
for suitable constants $y,s_0,c\in\R$. Moreover,  if we write
$$x(s)=\left(\sqrt{1+y^2}\right) \cosh(s+s_0) ,$$
we get by {\bf (C.4)}, \eqref{laache} and $\lambda u^2=1$ that $$ u=1/x(s), \quad\quad
\lambda=x(s)^2 .$$ Finally, observe that the vertical projection $N:\Sigma\flecha \H^2\subset \L^3$
of a mean curvature one half surface $\psi=(N,h):\Sigma\flecha \H^2\times \R$ verifies the
conditions
 \begin{equation*}
 \esiz N,G\esde = \frac{-1}{u},\hspace{.5cm} \esiz N,G_z\esde = \frac{1}{2}\left(h_z- \frac{2Q h_{\bar{z}}}{\tau}\right), \hspace{.5cm} \esiz N,N\esde =-1,
 \end{equation*}
where the second formula comes from \eqref{hisa}. Therefore we conclude that the vertical
projection of $\psi=(N,h)$ is given by $$\left\{\def\arraystretch{1.2} \begin{array}{lll}
N_0 & = & \displaystyle x(s)\cosh t + y \sinh t, \\
N_1 & = & x(s) \sinh t +y \cosh t, \\
N_2 & = &\pm\sqrt{ x(s)^2 - y^2 -1}.\end{array}\right. $$
As a result, the surface is invariant by hyperbolic screw motions, and has a simple explicit parametrization. This type of surfaces has been obtained in a totally different way by Sa Earp \cite{SaE}. They are entire vertical graphs over the whole horizontal factor $\H^2$ of $\H^2\times \R$. In particular, they are complete, embedded and stable.
\end{proof}

\subsection*{Parallel surfaces}

Let $\psi=(N,h):\Sigma\flecha \H^2 \times \R$ denote a surface with regular vertical projection, for which $H=1/2$ holds with respect to its canonical orientation (given by $u>0$, where $u$ is its angle function). Let $G:\Sigma\flecha \H^2$ denote its hyperbolic Gauss map, with Weierstrass data $\{-Q,\tau_0 =\landa u^2\}$. Here $\landa$ is the conformal factor of the metric of $\psi$, and $Q$ is its Abresch-Rosenberg differential.

Assume that $Q$ never vanishes. Then, by our discussion just after Theorem \ref{th:main} we know
that there exists another surface (actually a $2$-parameter family of them in the generic case)
$\psi^{\sharp} :\Sigma\flecha \H^2\times \R$ constructed via $G$ and the Weierstrass data
$\{-Q,\tau^{\sharp} = 16|Q|^2 /\tau_0\}$. This surface has regular vertical projection, it has
$H=1/2$ with respect to its canonical orientation, and its hyperbolic Gauss map is $G$.

It is then an interesting problem to establish if there is some explicit geometric relation between $\psi$ and $\psi^{\sharp}$. This is actually the case. The following theorem proves that $\psi$ and $\psi^{\sharp}$ are \emph{parallel} in a certain sense. This reproduces to some point the situation for parallel CMC $1/2$ surfaces in $\L^3$ with the same Gauss map that we exposed in Section 2.

\begin{teo}\label{teopa}
Let $\psi:\Sigma\flecha \H^2\times \R$ be a mean curvature one half surface with hyperbolic Gauss map $G:\Sigma\flecha \H^2$, and assume that its Abresch-Rosenberg differential $Q$ never vanishes. If $u:\Sigma\flecha (0,+\8)$ denotes its angle function, then
 \begin{equation}\label{paral}
 \psi^{\sharp} = -\psi + \frac{2}{u} (G,1) \hspace{1cm} \left( = -\psi + \frac{2}{u^2} (\eta + N)\right)
 \end{equation}
is a regular surface in $\H^2\times \R$ with regular vertical projection, for which $H=1/2$ holds for its canonical orientation. Moreover, its hyperbolic Gauss map is $G$, its angle function is $u$, and its conformal metric factor is
 \begin{equation}\label{landapar}
 \landa^{\sharp} = \frac{16|Q|^2}{\landa u^4}.
 \end{equation}
\end{teo}
In particular, this theorem tells that $\psi^{\sharp}$ as in \eqref{paral}   is the mean curvature
one half surface constructed according to Theorem \ref{th:main} by means of $G$ and the Weierstrass
data $\{-Q, \tau^{\sharp} = 16|Q|^2 /\tau_0\}$, with the initial condition (see the formula
\eqref{otracte} below) $$(h^{\sharp})_z (z_0)= \frac{4Q (z_0)}{\tau_0 (z_0)} \ \overline{h_0}.$$

\begin{proof}
Write $\psi=(N,h):\Sigma\flecha \H^2\times \R$, and let $\xi:\Sigma\flecha \N^3 $ be the map given by \eqref{sita}, where $\eta :\Sigma\flecha \S_1^3$ is the unit normal of $\psi$. Also recall that $\xi =(G,1)$. The following metric relations will be used repeatedly in what follows:
 \begin{equation}\label{parauno}
 \esiz \xi,\xi \esde =  \esiz \xi_z,\xi \esde =  \esiz \psi_z,\xi \esde =0, \hspace{0.5cm}  \esiz G,N \esde =  \esiz \xi,N \esde = -\frac{1}{u}.
 \end{equation}
Let us denote $\psi^{\sharp} = (N^{\sharp},h^{\sharp})$, where by \eqref{paral}
 \begin{equation}\label{parados}
 N^{\sharp} = -N + \frac{2}{u} G, \hspace{1cm} h^{\sharp} = -h + \frac{2}{u}.
 \end{equation}
By \eqref{parauno} and the fact that $G,N$ take values in $\H^2$, it follows directly that $\esiz N^{\sharp},N^{\sharp}\esde = -1$. Moreover, $\esiz G,N^{\sharp} \esde = -1 /u <0$, what ensures that $N^{\sharp}$ takes its values in $\H^2$. So, $\psi^{\sharp}$ is indeed a surface in $\H^2\times \R$, possibly with singular points.

Now, observe that by \eqref{paral} we have
 \begin{equation}\label{paratres}
 (\psi^{\sharp} )_z = -\psi_z + \left(\frac{2}{u}\right)_z \xi + \frac{2}{u}\,  \xi_z.
  \end{equation}
Thus, using \eqref{temita}, \eqref{eq:xiz}, \eqref{parauno}, \eqref{paratres} and $\esiz
\psi_z,\psi_z\esde =0$ we have
 \begin{equation}
 \esiz (\psi^{\sharp} )_z,(\psi^{\sharp} )_z\esde = \frac{4}{u^2} \esiz \xi_z,\xi_z\esde - \frac{4}{u} \esiz \psi_z,\xi_z\esde = \frac{4}{u^2} (-Q +Q)=0.
 \end{equation}
In other words, $z$ is a conformal parameter for $\psi^{\sharp}$, i.e. the surfaces $\psi$ and $\psi^{\sharp}$ have the same conformal structure.

Observe now that $\esiz (\psi^{\sharp} )_z,\xi\esde =0$ follows directly from \eqref{parauno} and
\eqref{paratres}. So, $\xi$ is a null vector in $\L^4$, whose last coordinate equals $1$, and that
is normal to $\psi^{\sharp}$. Therefore, if $\eta^{\sharp}:\Sigma\flecha \S_1^3$ is the unit normal
of $\psi^{\sharp}$ and $u^{\sharp}$ is its last coordinate (i.e. the angle function of
$\psi^{\sharp}$), we must have
 \begin{equation*}
 \xi = \frac{1}{u^{\sharp}} (\eta^{\sharp} + \ep N^{\sharp} ), \hspace{1cm} \ep = \pm 1.
 \end{equation*}
Consequently we get $$ \esiz N^{\sharp},\xi\esde = \esiz N^{\sharp},\frac{1}{u^{\sharp}}
(\eta^{\sharp} + \ep N^{\sharp} )\esde = -\frac{\ep}{u^{\sharp}}.$$ On the other hand, by
\eqref{parados} we see that $$ \esiz N^{\sharp},\xi\esde = \esiz -N +\frac{2}{u} \, G,G\esde =
-\frac{1}{u}.$$ Thus $u^{\sharp} =\ep u$. In particular, $\psi^{\sharp}$ has regular vertical
projection. So, if we endow $\psi^{\sharp}$ with its canonical orientation, given by $u^{\sharp}
>0$, we have
 \begin{equation}\label{paracinco}
 u^{\sharp} =u \hspace{0.6cm} \text{ and } \hspace{0.6cm} \xi = \frac{1}{u^{\sharp}} (\eta^{\sharp} + N^{\sharp} ).
 \end{equation}
In other words, $G:\Sigma\flecha \H^2$ is the hyperbolic Gauss map of $\psi^{\sharp}$.

The conformal factor $\landa^{\sharp}$ of $\psi^{\sharp}$ is obtained directly from \eqref{parauno}, \eqref{paratres} and \eqref{eq:xiz}, as follows:

\begin{equation}\def\arraystretch{2.2}\begin{array}{lll}
\landa^{\sharp} & = & \displaystyle 2 \esiz (\psi^{\sharp} )_z, (\psi^{\sharp} )_{\bar{z}}\esde =
\landa -\frac{8}{u}\esiz \psi_z,\xi_{\bar{z}}\esde + \frac{8}{u^2} \esiz \xi_z,\xi_{\bar{z}}\esde \\
& =& \landa \displaystyle -\frac{8}{u} \frac{\landa u}{4} + \frac{8}{u^2} \left(\frac{\landa
u^2}{8} + \frac{2|Q|^2}{\landa u^2}\right) = \frac{16|Q|^2}{\landa u^4}.
\end{array}\end{equation}
That is, \eqref{landapar} holds. Moreover, as $Q$ never vanishes, $\psi^{\sharp}$ is regular.

Finally, we just need to check that $\psi^{\sharp}$ has constant mean curvature $H^{\sharp} =1/2$.
For that, let us observe first of all that by ${\bf (C.3)}$ with $H=1/2$ and ${\bf (C.4)}$, and
using that $Q =p +h_z^2$, we have $$\left(\frac{1}{u}\right)_z = \frac{h_z}{2} + \frac{2 Q}{\landa
u^2} h_{\bar{z}}.$$ So, since by \eqref{parados} it holds $h^{\sharp} = -h + 2/u$, we conclude that

 \begin{equation}\label{otracte}
(h^{\sharp})_z = \frac{ 4Q}{\landa u^2} h_{\bar{z}}.
 \end{equation}
If we differentiate this expression, then by ${\bf (C.1)}$, ${\bf (C.3)}$, ${\bf (C.4)}$ and
$Q=p+h_z^2$ we get
\begin{equation}\label{paraseis}
(h^{\sharp})_{z\bar{z}} = \frac{4 |Q|^2}{\landa u^3}.
\end{equation}
By ${\bf (C.2)}$, the mean curvature $H^{\sharp}$ of $\psi^{\sharp}$ is $$H^{\sharp} = \frac{2 (h^{\sharp})_{z\bar{z}}}{\landa^{\sharp} u^{\sharp}}.$$ At last, using \eqref{landapar}, \eqref{paracinco} and \eqref{paraseis} we arrive at $H^{\sharp} = 1/2$. This ends up the proof.
\end{proof}

%Let us make the following remark regarding this construction of parallel surfaces. In $\L^3$ the only CMC $1/2$ surface that is \emph{self-parallel} (i.e. it agrees with its parallel CMC $1/2$ surface) is the hyperbolic cylinder. This follows directly from \eqref{yanise}. It is also the only CMC $1/2$ surface with everywhere singular Gauss map. In our situation of $\H^2\times \R$, the same happens for the the standard hyperbolic screw motion examples of Sa Earp. Indeed, it follows from \eqref{yanise} and Proposition \ref{sae} that these examples are the only mean curvature one half surfaces that are self-parallel, i.e. $\psi^{\sharp} = \psi$. In addition, we saw in Proposition \ref{sae} that they are also the only surfaces with singular hyperbolic Gauss map. So, it seems that these screw motion surfaces with $H=1/2$ in $\H^2\times \R$ may play to some extent the role of hyperbolic cylinders in CMC surface theory in $\L^3$.

\subsection*{Surfaces with conformal hyperbolic Gauss map}

We will construct next all the mean curvature one half surfaces in $\H^2\times \R$ with regular vertical projection whose hyperbolic Gauss map is conformal. A way to do this relies on the fact that the system \eqref{eq:A} can be explicitly integrated when $Q\equiv 0$. However, we will instead use a shorter argument based on the \emph{parallel surfaces} construction of the previous subsection. Let us anyway indicate that the CMC surfaces with $Q=0$ in $\H^2\times \R$ and $\S^2\times \R$ were classified in \cite{AbRo}.

Let $\psi =(N,h):\Sigma\flecha \H^2\times \R$ be a mean curvature one half surface with regular vertical projection. Assume, besides, that its hyperbolic Gauss map $G:\Sigma\flecha \H^2$ is conformal, i.e. $\esiz G_z,G_z\esde =0$. Thus, the Abresch-Rosenberg differential $Q dz^2$ of $\psi$ vanishes identically. This indicates by \eqref{landapar} that its parallel surface $\psi^{\sharp}:\Sigma\flecha \H^2\times \R$ given by \eqref{paral} verifies $\esiz d\psi^{\sharp},d\psi^{\sharp}\esde =0$. In other words, $\psi^{\sharp}$ is constant. It is also clear that this property characterizes the surfaces with $H=1/2$ and conformal hyperbolic Gauss map.

Write now $\psi^{\sharp} =(a,b)\in \H^2\times \R$. We plan to recover $\psi$ in terms of $a$ and
$G$.  By \eqref{paral} we obtain $$ \psi = -(a,b) + \frac{2}{u} (G,1).$$ Thus, by applying if
necessary a vertical translation to $\psi$ we see that $h=2/u$, and hence
 \begin{equation}\label{conf1}
 \psi =(N,h)= -(a,0) + h (G,1).
 \end{equation}
At last, as $G\in \H^2$ and $a\in \H^2$, the condition $\esiz N,N\esde =-1$ provides
 \begin{equation}\label{conf2}
 h=-2\esiz a,G\esde.
 \end{equation}
Putting \eqref{conf1} and \eqref{conf2} together we obtain that $\psi$ is expressed in terms of the conformal map $G:\Sigma\flecha \H^2$ and a constant $a\in \H^2$ as
 \begin{equation}\label{conf3}
 \psi = (-a,0) - 2\esiz a,G\esde (G,1):\Sigma\flecha \H^2\times \R\subset \L^4 .
 \end{equation}
Conversely, if $G:\Sigma\flecha \H^2$ is conformal and regular, and $a\in \H^2$, then the map $\psi$ given by \eqref{conf3} is a mean curvature one half surface in $\H^2\times \R$ with regular vertical projection, having $G$ as its hyperbolic Gauss map. We omit the proof, as it is a direct computation using the ideas in the proof of Theorem \ref{teopa}.

\subsection*{A reflection principle}

We show next that mean curvature one half surfaces in $\H^2\times \R$ admit a Schwarz reflection principle. It roughly states that a surface with $H=1/2$ in $\H^2\times \R$ that meets orthogonally a totally geodesic vertical plane $\mathcal{P}$ of $\H^2\times \R$ can be analytically extended by reflection across $\mathcal{P}$ as a surface with $H=1/2$.

Let $\Omega\subseteq \C$ be a complex symmetric domain, i.e. $\Omega =\Omega^* :=\{\bar{z}: z\in \Omega\}$. Define $\Omega^+ = \Omega \cap \{{\rm Im} (z) >0\}$, $\Omega^- = \Omega \cap \{{\rm Im} (z) <0\}$, and $I= \Omega\cap \R$.

Let $\mathcal{P}$ denote a totally geodesic vertical plane of $\H^2\times \R$, i.e. a right vertical cylinder over a geodesic of $\H^2$. Finally, define $\sigma$ as the isometric reflection in $\H^2\times \R$ with respect to $\mathcal{P}$.

Then we have

\begin{teo}
Let $\psi:\Omega^+\cup I \subset \C\flecha \H^2\times \R$ be a conformal $C^2$  immersion with $H=1/2$ and regular vertical projection. Assume that $\psi$ maps $I$ into $\mathcal{P}$ and meets $\mathcal{P}$ orthogonally, i.e. its unit normal $\eta :\Omega^+ \cup I\flecha \S_1^3$ is tangent to $\mathcal{P}$ along $I$.

Then the map $\psi:\Omega \subset\C\flecha \H^2\times \R$ given by
 \begin{equation*}
 \psi (z) =\left\{\def\arraystretch{1.4}\begin{array}{lll} \psi (z) & \text { if } & z\in \Omega^+\cup I , \\
\sigma( \psi (\bar{z})) & \text { if } & z\in \Omega^-
 \end{array} \right.
 \end{equation*}
is a conformal immersion with $H=1/2$ which extends $\psi$ symmetrically across $\mathcal{P}$.
\end{teo}
\begin{proof}
Up to a rigid motion, we can take $\mathcal{P}$ to be $\mathcal{P}= (\H^2 \times \R) \cap \{x_2=0\} \subset \L^4$, and thus we have $\sigma (x_0,x_1,x_2,x_3)= (x_0,x_1,-x_2,x_3)$. Then, the hypothesis on $\psi$ indicates that the hyperbolic Gauss map $$G=(G_0,G_1,G_2):\Omega^+ \cup I \flecha \H^2 \subset \L^3$$ of $\psi$ is continuous and verifies that $G_2(s,0)=0$ for all $s\in I$. In the same way, we see that
 \begin{equation}\label{refl1}
 \frac{\parc \psi}{\parc t} (s,0)= \left( 0,0,\frac{\parc \psi_2}{\parc t} (s,0), 0\right)
 \end{equation}
for all $s\in I$. Now let $\xi =(G_0,G_1,G_2,1)$. As $G_2 (s,0)=0$, the third coordinate of \eqref{eq:xiz} at $s\in I$ provides $$ \frac{\parc G_2}{\parc t} = \left(\frac{u}{2} - \frac{2 Q}{\landa u}\right) \frac{\parc \psi_2}{\parc t}.$$ In particular $Q(s,0)\in \R$ for all $s\in I$. Using this fact we infer that the imaginary part of the first two coordinates of \eqref{eq:xiz} at $s\in I$ turn into $$\frac{\parc G_0}{\parc t} (s,0)=0, \hspace{1cm} \frac{\parc G_1}{\parc t} (s,0)=0.$$ Therefore, we have seen that the harmonic map $G:\Omega^+\cup I\flecha \H^2$ verifies that $G(I)$ is a part of the geodesic $\gamma=\H^2\cap \{x_2=0\}$ of $\H^2$, and the tangential component of the normal derivative $(\parc G/\parc t) (s,0)$ vanishes. Hence, by the Schwarz reflection principle for harmonic maps (see \cite{Woo}) we conclude that $G$ can be harmonically extended to $\Omega$ by means of
\begin{equation}\label{refl3}
G(z)= \left\{\def\arraystretch{1.4}\begin{array}{lll} G (z) & \text { if } & z\in \Omega^+ \cup I, \\
G(\bar{z})^* & \text { if } & z\in \Omega^- ,
 \end{array} \right.
 \end{equation}
where here $P^*$ denotes the reflection of $P\in \H^2$ across $\gamma \subset\H^2$.

As a result, by Theorem \ref{th:main}, we can extend the surface $\psi$ to a conformal immersion $\psi:\Omega \flecha \H^2\times \R$ with $H=1/2$. Moreover, let $\Psi:\L^3\flecha \L^3$ be the symmetry in $\L^3$ that extends the isometric reflection in $\H^2$ across $\gamma$. Then $\Psi (G(\bar{z})) = G(z)$, and so $\sigma (\psi (\bar{z})) =\psi (z)$, where $\sigma= (\Psi,{\rm Id})$ is the symmetry of $\H^2 \times \R$ with respect to $\mathcal{P}$. This completes the proof.

\end{proof}
\subsection*{Minimal surfaces in the Heisenberg group}

In \cite{Dan} the classical Lawson correspondence between CMC surfaces in space forms was extended to CMC surfaces in other $3$-dimensional homogeneous spaces. This generalized Lawson correspondence shows, in particular, the existence of a bijective isometric correspondence between simply connected surfaces with $H=1/2$ in $\H^2\times \R$ and simply connected minimal surfaces in the $3$-dimensional Heisenberg group ${\rm Nil}_3$.

In this final part of the section we exploit this fact in order to obtain an existence result for complete simply connected minimal surfaces in ${\rm Nil}_3.$ First we will explain briefly this correspondence.

The space ${\rm Nil}_3$ can be regarded as the Lie group
$${\rm Nil}_3=\left\{ \begin{pmatrix}
1&x_1&\frac{1}{2}x_1x_2+x_3\\
0&1&x_2\\
0&0&1\end{pmatrix}\;:\; (x_1,x_2,x_3)\in\R^3
\right\}, $$
endowed with the left invariant metric
$$ds_0^2=dx_1^2+dx_2^2+\left(\frac{1}{2} (x_2dx_1-x_1dx_2) + dx_3\right)^2.$$
We denote by  $\chi=\partial_{x_3}$
the Killing vector field corresponding to the {\em vertical translations} in ${\rm Nil}_3.$\\

Let $\psi=(N,h):\Sigma\to\H^2\times\R$ be a mean curvature one half immersion from an open simply connected Riemann surface $\Sigma$, with first fundamental form $ds^2=\lambda|dz|^2$  and Hopf differential $p\,dz^2.$
Let $\chi_0$ be the Killing field corresponding to the vertical translations in $\H^2\times\R$, and denote by $\eta=(\hat{N},u)$ the unit normal of $\psi.$ It easy to check that
$$\chi_0=\frac{2}{\lambda}(h_{\bar z}\psi_z + h_{z}\psi_{\bar z}) +  u \eta. $$

Then, by the generalized Lawson correspondence in \cite{Dan}, there exists a (unique up to rigid motions) conformal minimal immersion $\widehat{\psi}:\Sigma\to{\rm Nil}_3$ with first fundamental form $ds^2 =\landa |dz|^2,$ with Hopf differential $\widehat{p}\,dz^2=i p\, dz^2$ and satisfying
$$\chi=\frac{2i}{\lambda}(h_{\bar z}\widehat{\psi}_z + h_{z}\widehat{\psi}_{\bar z}) +  u \widehat{\eta}, $$
where here $\widehat{\eta}$ is the unit normal of $\widehat{\psi}.$ Such a  pair of immersions $(\psi,\widehat{\psi})$ are called {\em sister surfaces} \cite{Dan}.

Recently, Abresch \cite{Abr} (see also \cite{BeTa,FeMi}) has announced the existence of a holomorphic quadratic differential for CMC surfaces in the Heisenberg group ${\rm Nil_3}$. Using the above notation, and in the case of minimal surfaces, this differential can be written as
$$\widehat{Q}= i\, \widehat{p}\,dz^2  + \langle \chi , \widehat{\psi}_z\rangle^2 dz^2 . $$

Thus, the relation between this differential and the Abresch-Rosenberg differential $Q$ for its sister surface of $H=1/2$ in $\H^2\times\R$ is given by $\widehat{Q}=-Q.$ Keeping this in mind we obtain the following corollary to Theorem \ref{exicom}.

\begin{cor}
Any holomorphic quadratic differential on an open simply connected Riemann surface $\Sigma$ is the Abresch differential of some complete minimal surface in ${\rm Nil}_3$. Furthermore, the space of congruency classes of complete minimal surfaces in ${\rm Nil}_3$ with the same Abresch differential is generically infinite.
\end{cor}

\section{Complete minimal graphs in $\H^2\times \R$}

Let $\psi =(N,h):\Sigma\flecha \H^2\times \R$ be an orientable minimal surface, and consider $\Sigma$ as a Riemann surface with the conformal structure given by its metric. By {\bf (C.2)} we see that the height function $h$ is harmonic, and thereby $\omega := h_z dz$ is a globally defined holomorphic $1$-form on $\Sigma$, where $z$ denotes an arbitrary complex coordinate on $\Sigma$. We call $\omega$ the \emph{canonical $1$-form} of the minimal surface $\psi$.

Consider also the unit normal $\eta:\Sigma\flecha \S_1^3 \subset \L^4$ of $\psi$. We will define as usual the \emph{angle function} $u:\Sigma\flecha [-1,1]$ as the last coordinate of $\eta$, i.e. $u=\esiz \eta,e\esde$ where $e=(0,0,0,1)$. Observe that $u$ and $\omega$ are closely related by {\bf (C.4)}.

If we write $\esiz d\psi,d\psi \esde =\landa |dz|^2$ for a positive smooth function $\landa$ on $\Sigma$, then the metric of the vertical projection $N$ of $\psi$ is
 \begin{equation}\label{enemu}
 \esiz dN,dN\esde = -h_z^2 dz^2 + \mu |dz|^2 - h_{\bar{z}}^2 d\bar{z}^2, \hspace{0.5cm} \mu = \landa - 2|h_z|^2.
  \end{equation}
In particular, $N:\Sigma\flecha \H^2$ is harmonic. Observe also that $dN\neq 0$ at every point, since $\landa$ is positive.

Conversely, if $N:\Sigma\flecha \H^2$ is a harmonic map with $dN\neq 0$ everywhere, and whose Hopf differential is of the form $Q_0 dz^2=-\omega^2$ for some holomorphic $1$-form $\omega$ on $\Sigma$, then
\begin{equation}\label{repmin}
\psi := \left( N,2\Re \int \omega \right) :\Sigma\flecha \H^2\times \R
\end{equation}
is a conformal minimal surface, provided that the above integral has no real periods.

With this, we have:

\begin{teo}\label{supor}
Let $\mathcal{A}$ denote the space of complete orientable minimal surfaces in $\H^2\times \R$ whose angle function omits some value in $(-1,1)$. Then:
 \begin{enumerate}
 \item[{\rm (i)}]
The elements of $\mathcal{A}$ are exactly the complete minimal vertical graphs, and the right vertical cylinders over some geodesic of $\H^2$.
 \item[{\rm (ii)}]
All the elements of $\mathcal{A}$ are simply connected.
 \item[{\rm (iii)}]
Two elements $\psi_1,\psi_2 \in \mathcal{A}$ are congruent if and only if their respective canonical $1$-forms $\omega_1,\omega_2$ verify $\omega_1 = \pm \omega_2$.
 \item[{\rm (iv)}]
Every holomorphic $1$-form on $\D$, and every non-zero holomorphic $1$-form on $\C$, can be realized as the canonical $1$-form of some element of $\mathcal{A}$.
 \end{enumerate}
\end{teo}
 \begin{proof}
Observe first of all that every complete minimal vertical graph in $\H^2\times \R$ lies in $\mathcal{A}$, since its angle function $u$ does not change sign. Also note that the right vertical cylinders over geodesics of $\H^2$ have vanishing angle function, and thus are also in $\mathcal{A}$.

Let $\psi:\Sigma\flecha \H^2\times \R$ be an element of $\mathcal{A}$, and consider the function $$\sigma:= \landa (1+u)^2 :\Sigma\flecha [0,+\8).$$ By the hypothesis on the angle function $u$, and composing if necessary with an inverse rigid motion in $\H^2\times \R$, we may assume that $u\geq c>-1$ for some $c\in (-1,1]$. Thus $\sigma$ is positive, and we have
 \begin{equation}\label{acotmin}
(1+c)^2 \landa \leq \sigma \leq 4 \landa.
 \end{equation}
Hence, the completeness of the metric $\sigma |dz|^2$ is equivalent to the completeness of the surface $\psi$.

Moreover, using ${\bf (C.1)}$ to ${\bf (C.4)}$ and Codazzi in \eqref{GCR} with $H=0$, we see that $$\def\arraystretch{2.3}\begin{array}{lll}
(\log (1+u))_{z\bar{z}} & = & \displaystyle\frac{u_{z\bar{z}}}{1+u} -\frac{u_z u_{\bar{z}}}{(1+u)^2} =  \frac{u |h_z|^2}{1+u} - \frac{2 |p|^2 u}{\landa (1+u)} -\frac{|p|^2 (1-u^2)}{\landa (1+u)^2} \\ & = & \displaystyle \frac{\landa u (1-u)}{4} - \frac{ |p|^2}{\landa}. \end{array}$$ Also observe that the Gauss equation in \eqref{GCR} for $H=0$ together with ${\bf (C.4)}$ imply $$ (\log \landa)_{z\bar{z}} = \frac{2|p|^2}{\landa} + \frac{\landa u^2 }{2}.$$ Therefore,
 \begin{equation*}%\label{sima1}
 (\log \sigma)_{z\bar{z}} = \frac{\landa u}{2}.
 \end{equation*}
With this, and by {\bf (C.4)} we have $$(\log \sigma)_{z\bar{z}} = \frac{\sigma}{8} - \frac{2|h_z|^4}{\sigma} \hspace{1cm} \left(  =  \frac{\landa u}{2}  \right) .$$ In this way, if we denote by $\tilde{\Sigma}$ the universal cover of $\Sigma$, and $\pi:\tilde{\Sigma}\to\Sigma$ is the corresponding covering map,
it follows that $\sigma\circ\pi$ is a solution to \eqref{eq:gauss}, the Gauss equation for $H=1/2$ surfaces in $\L^3$. So, there exists a unique (up to rigid motions in $\L^3$) spacelike surface $f:\tilde{\Sigma}\flecha \L^3$ with $H=1/2$, whose metric is $\esiz df,df\esde = (\sigma\circ\pi)|d\zeta|^2$ (hence, complete), and whose Hopf differential is
$Q_0 d\zeta^2 := -(\pi^*\omega)^2.$ Here $\zeta=z\circ\pi$ is a conformal parameter on $\tilde{\Sigma}$ and $\pi^*\omega$ denotes the pullback of the 1-form $\omega$ by the covering $\pi.$  Furthermore, let $\nu:\tilde{\Sigma}\flecha \H^2\cup \H^2_{-}\subset \L^3$ denote the Gauss map of $f$, and let $P$ be the rigid motion of $\L^3$ so that $P={\rm Id}$ if $\nu (\tilde{\Sigma})\in \H^2$ and $P(x_0,x_1,x_2)=(-x_0,x_1,x_2)$ if $\nu (\tilde{\Sigma})\in \H^2_{-}$. Here $\H^2_{-} = \{(-x_0,x_1,x_2) \in \L^3 : (x_0,x_1,x_2)\in \H^2\subset \L^3\}$. Then, if we define $g= P\circ \nu:\tilde{\Sigma}\flecha \H^2$, by \eqref{metG} and \eqref{enemu} we have $$\langle g_{\zeta},g_{\zeta}\rangle =
Q_0 = \langle (N\circ\pi)_{\zeta},(N\circ\pi)_{\zeta}\rangle$$
and
$$\langle g_{\zeta},g_{\bar \zeta}\rangle = \frac{\sigma\circ\pi}{8}+\frac{2|Q_0|^2}{\sigma\circ\pi}=
\left(\frac{\lambda(1+u^2)}{4}\right)\circ\pi = \frac{\mu}{2}\circ\pi  =\langle (N\circ\pi)_{\zeta},(N\circ\pi)_{\bar\zeta}\rangle.$$
So, there exists an isometry $\Phi$ of $\H^2$ so that $\Phi \circ g = N\circ \pi$, i.e. $\Psi \circ \nu = N \circ \pi$ where $\Psi = \Phi \circ P$.

Recall now (see \cite{ChYa,ChTr,Wan}) that, by the completeness of $f$, the Gauss map
$\nu:\tilde{\Sigma}\flecha \H^2$ is a global diffeomorphism from $\tilde{\Sigma}$ onto
$\nu(\tilde{\Sigma})$, except if the surface is a hyperbolic cylinder in $\L^3$, in which case
$\nu$ is a piece of a geodesic of $\H^2$. In the first situation, by $\Psi\circ \nu =N\circ \pi$,
we necessarily conclude that $\psi$ is a vertical graph in $\H^2\times \R$, and that $\pi$ is
one-to-one, that is, $\Sigma$ is simply connected. In the second one, again by $\Psi\circ \nu
=N\circ \pi$, $\psi$ must be a right vertical cylinder over a geodesic of $\H^2$. This proves (i)
and (ii). These results let us assume in the remaining part of the proof that $\Sigma\equiv
\tilde{\Sigma}$ and $\pi = {\rm Id}$.

To prove (iii), recall first that by \cite{Wan} and \cite{WaAu}, a complete $H=1/2$ surface in $\L^3$ is uniquely determined by its Hopf differential. So, using $\Psi\circ \nu =N$, we get that two elements of $\mathcal{A}$ are congruent if and only if they produce following the above process two congruent complete $H=1/2$ surfaces in $\L^3$. Recalling finally that the canonical $1$-form $\omega$ of $\psi$ and the Hopf differential $Q_0$ of $f$ are related by $Q_0 =-\omega^2$, we obtain (iii).

Finally, to prove the existence part (iv), let us start with a holomorphic $1$-form $\omega$ on $\Sigma = \C$ or $\D$, in the conditions of the theorem. By \cite{Wan} and \cite{WaAu} there is a unique (up to rigid motions) complete spacelike $H=1/2$ surface $f:\Sigma\flecha \L^3$ whose Hopf differential is $Q_0:=-\omega^2$.

Let $\nu:\Sigma\flecha \H^2 \cup \H^2_{-}$ and $\tau_0$ denote, respectively, the Gauss map and the conformal factor of the metric of $f$. Define now $N=P\circ \nu :\Sigma\flecha \H^2$, where $P$ is the rigid motion of $\L^3$ defined as above. Then we can construct a minimal surface $\psi=(N,h):\Sigma\flecha \H^2\times \R$ by means of the representation formula \eqref{repmin}. By its construction, $\psi$ is a vertical graph, or a right vertical cylinder over a geodesic of $\H^2$, and its canonical $1$-form is $\omega$. Now, putting together \eqref{enemu} and \eqref{metG} we infer that $$\landa = \mu + 2|h_z|^2 = \frac{\tau_0}{4} + \frac{4 |h_z|^4}{\tau_0} + 2|h_z|^2 = \frac{(\tau_0 + 4|h_z|^2)^2}{4\tau_0} \geq \frac{\tau_0}{4}.$$ Hence $\psi$ is regular and complete, and so $\psi\in \mathcal{A}$. This ends up the proof.
 \end{proof}

The classification of the complete minimal vertical graphs in $\H^2\times \R$ follows readily from Theorem \ref{supor}.

\begin{cor}
Consider the following spaces:

$$ \def\arraystretch{1.1}\begin{array}{lll}
\mathcal{G} &= & \{ \text{congruency classes of complete minimal vertical graphs in $\H^2\times \R$}\}.\end{array}$$

$$ \begin{array}{lll}
\Theta &= & \{ \text{holomorphic $1$-forms on $\C $ which are not of the type $\omega = c\, dz$, $c\in \C$,} \\ & & \text{ and holomorphic $1$-forms on $\D$}
 \}.\end{array}$$
Let finally $\Theta /\Z_2$ denote the quotient of $\Theta$ obtained by identifying elements of $\Theta$ differing only by a $\pm $ sign.

Then the map assigning to each element $g\in \mathcal{G}$ the class in $\Theta /\Z_2$ of the canonical $1$-form of some representant $\psi:\Sigma\flecha \H^2\times \R$ of $g$ is a well defined bijective correspondence between $\mathcal{G}$ and $\Theta /\Z_2$.
\end{cor}
\begin{proof}
Observe first of all that the conformal structure is invariant by rigid motions, and that by Theorem \ref{supor}, all complete minimal vertical graphs in $\H^2\times \R$ are conformally equivalent to $\C$ or $\D$.

Also, the canonical $1$-form of any complete minimal vertical graph in $\H^2\times \R$ lies in $\Theta$. Indeed, by Theorem \ref{supor}, an element of $\mathcal{A}$ fails to be a vertical graph if and only if it is a right vertical cylinder over a geodesic, if and only if its associated complete $H=1/2$ surface in $\L^3$ is a hyperbolic cylinder, if and only if its conformal structure is parabolic (hence $\Sigma\equiv \C$), and $Q_0=a \,dz^2$ for some $a\in \C\setminus\{0\}$ (this is by the completeness of the $H=1/2$ associated surface), if and only if $\Sigma\equiv \C$ and $\omega = c \, dz$ for some $c\in \C\setminus\{0\}$.

With this, by part (iii) in Theorem \ref{supor}, the correspondence described by the present theorem is well defined, and it also follows directly from Theorem \ref{supor} that it is a bijection.
\end{proof}

The previous results let us also draw the following conclusion: \emph{a harmonic map from an open Riemann surface $\Sigma$ into $\H^2$ is the vertical projection of a complete minimal vertical graph in $\H^2\times \R$ if and only if: (1) its Hopf differential is of the form $Q dz^2 = -\omega^2$ for some holomorphic $1$-form $\omega$ on $\Sigma$, and (2) it is the Gauss map of a complete spacelike CMC surface in $\L^3$ different from a hyperbolic cylinder}.

Therefore, several theorems regarding the Gauss map image of complete CMC surfaces in $\L^3$ can be translated into our context. For example, by the results in \cite{Aiy,CSZ,Xin,XiYe} we have:
 \begin{cor}
Let $M^2\subset \H^2\times \R$ be a complete minimal vertical graph over a domain $\Omega \subset \H^2$. Then
 \begin{enumerate}
 \item
$\Omega$ cannot be bounded.
 \item
$\Omega$ cannot lie inside a tubular neighborhood of a geodesic in $\H^2$.
 \item
$\Omega$ cannot lie inside a horoball.
 \end{enumerate}
 \end{cor}

Furthermore, using \cite{HTTW} we obtain the following consequence:
 \begin{cor}
Let $M^2\subset \H^2\times \R$ denote a complete minimal vertical graph over a domain $\Omega \subset \H^2$, and suppose that $M^2$ has parabolic conformal type. Then $\Omega$ is an ideal geodesic polygon with $n$ vertices in $\H^2$ if and only if $n$ is even and the canonical $1$-form of $M^2$ is $\omega = p(z) dz$, where $p(z)$ is a polynomial of degree $n/2$.
 \end{cor}

\def\refname{References}

\end{document}